\documentclass[12pt]{amsart}
\usepackage{mschoenbc}
\usepackage{relsize}

\begin {document}
\title{
 Exact Formulas for Invariants of Hilbert Schemes}

\author{Nate Gillman}
\address{Department of Mathematics \& Computer Science, Wesleyan University, Middletown, CT 06457, U.S.A.}
\email{ngillman@wesleyan.edu}

\author{Xavier Gonzalez}
\address{Mathematical Institute, University of Oxford, Oxford, England}
\email{{xavier.gonzalez@balliol.ox.ac.uk}}

\author{Matthew Schoenbauer}
\address{Department of Mathematics, University of Notre Dame, Notre Dame, IN 46556, U.S.A.}
\email{mschoenb@nd.edu}

\date {\today}
\maketitle


\begin{abstract}
  A theorem of G\"ottsche establishes a connection between cohomological invariants of a complex projective surface $S$ and corresponding invariants of the Hilbert scheme of $n$ points on $S$. 
  This relationship is encoded in certain infinite product $q$-series which are essentially modular forms. 
  Here we make use of the circle method to arrive at exact formulas for certain specializations of these $q$-series, yielding convergent series for the signature and Euler characteristic of these Hilbert schemes.
  We also analyze the asymptotic and distributional properties of the $q$-series' coefficients.
\end{abstract}

\section{Introduction and Statement of Results}

K3 surfaces are complex surfaces  characterized by their particularly simple Hodge structures and trivial holomorphic tangent bundles. 
These manifolds are of particular interest to physicists and mathematicians.
For physicists, they serve as useful examples of Calabi-Yau manifolds, which are a class of spaces central in string theory, and for mathematicians, they serve as an interesting yet sufficiently simple example in 4-manifold theory and complex differential geometry. 
These two roles came together in an unexpected way when in \cite{YZ} Yau and Zaslow conjectured (and later Beauville proved in \cite{Arn}) that the count of $n$-nodal curves on a K3 surface is equal to the Euler characteristic $\chi(\text{Hilb}^n(S))$ of Hilbert schemes of $n$ points on a K3 surface $S$. 

Yau and Zaslow's conjecture made use of a previous theorem of  G\"ottsche (see \cite[p.~37]{GotBook}) that provides an infinite product encoding the Euler characteristics $\chi(\mathrm{Hilb}^n(S))$ for Hilbert schemes of a K3 surface $S$. These Euler characteristics can be assembled in the form of the generating function 
\beq X_S(\tau) := \sum_{n =0}^\infty \chi(\mathrm{Hilb}^n(S)) q^n= \frac{q}{\Delta(\tau)} =\prod_{n=1}^\infty\frac{1}{(1-q^n)^{24}},\label{eq:Delta} \eeq
where $\Delta(\tau)$ is the modular discriminant and $q:=e^{2\pi i\tau}.$

G\"ottsche stated a more refined infinite product formula concerning Hodge numbers, which are cohomological invariants  that can be assembled into Betti numbers and the Euler characteristic. Recently, Manschot and Zapata Rolon \cite{Manschot} studied the asymptotic distribution of linear combinations of these Hodge numbers, which are realized as coefficients of Laurent polynomials called the $\chi_y$ genera for K3 surfaces (see (\ref{eq:chiy}) for a definition of $\chi_y$). The $\chi_y$ genera can also be assembled using the generating function \beq Y_S(y; \tau) := \sum_{m,n} b_S(m;n) y^m q^n := \sum_{n \geq 0} \chi_y(\mathrm{Hilb}^n(S)) q^n \label{eq:Y-series}. \eeq 
Akin to (\ref{eq:Delta}), G\"ottsche also provides an infinite product for $Y_S(y;\tau)$ (see Lemma \ref{lem:gottsche_prod}). 

Manschot and Zapata Rolon \cite{Manschot} found for  K3 surfaces $S$ that if $m$ is fixed, then as $n\rightarrow\infty$ we have
$$b_S(m;n) \sim \dfrac{\pi}{3 \sqrt{2}} n^{-\frac{29}{4}} \cdot \exp(4 \pi \sqrt{n}).$$
Moreover, since the asymptotics for $b_S(m;n)$ do not depend on $m,$ if we define  $$b^*_S(r; n) := \sum_{m \equiv r \Mod{2}} b_S(m;n), $$
it follows that $b^*_S(0;n)\sim b^*_S(1;n)$ as $n \to \infty.$
From a geometric perspective, $n$ corresponds to the length of the Hilbert scheme, and $m$ corresponds to a particular monomial in the $\chi_y$-genus of $\text{Hilb}^n(S)$. Grouping the contributions by the coefficients of these monomials in residue classes mod $2$, we obtain an equidistribution in the limit as $n$ goes to infinity.

G\"ottsche's formula for Hodge numbers $h^{s,t}(\text{Hilb}^n(S))$ of Hilbert schemes of $n$ points for K3 surfaces is a special case of his more general formula, which states that for any smooth  projective complex surface $S$, we have 
\beq \label{eq:zetadef} Z_S(x,y;\tau) 
:= \sum_{n\geq0}\chi_{\mathrm{Hodge}}(\mathrm{Hilb}^n(S))q^n =\prod_{n =1}^\infty \frac{\prod_{s +t \mathrm{\ odd}} (1- x^{s-1}y^{t-1}q^n)^{h^{s,t}}}{\prod_{s +t \mathrm{\ even}} (1- x^{s-1}y^{t-1}q^n)^{h^{s,t}}}. \eeq
For later use, we define $c_S(s,t;n)$ to be the coefficient of $x^sy^tq^n$ in the power series expansion of $Z_S(x,y;\tau) $.
\skipaline \noindent \tbf{Remark:} The \emph{Hodge polynomial} $\chi_{\mathrm{Hodge}}(\mathrm{Hilb}^n(S))$ is a Laurent polynomial in $\bZ[x,y,x^{-1}, y^{-1}],$ and we will sometimes write it as $\chi_{\mathrm{Hodge}}(\mathrm{Hilb}^n(S))(x,y)$ to make explicit the Hodge polynomial's dependence on $x$ and $y.$
If we specialize $x$ and $y$ to $\pm1$, then $\chi_{\mathrm{Hodge}}(\mathrm{Hilb}^n(S))(x,y)$ evaluates to different linear combinations of important topological invariants. 
See Section \ref{sec:gen_functions}  
for a more detailed discussion of $\chi_{\mathrm{Hodge}}(\mathrm{Hilb}^n(S))(x,y)$.
\skipaline

We seek exact formulas for  sequences assembled from the coefficients of (\ref{eq:zetadef}) for  a more general class of smooth projective complex surfaces. We consider the distribution of the coefficients  $c_S(s,t;n)$ over residue classes mod 2; namely, if we let
\beq \label{eq:c-special} c^*_{S}(r_1, r_2;n):= \sum_{\substack{{t \equiv r_1} \Mod{2}\\ {s \equiv r_2} \Mod{2} }} c_S(s,t;n),\eeq
we seek to determine the asymptotic properties of the sequences $b^*_S(r;n)$ and $c^*_{S}(r_1, r_2;n)$.
To this end, we consider the asymptotics of the $q$-series \beq \label{eq:bigC} C_S(r_1,r_2;\tau) := \sum_{n \geq 0} c^*_{S}(r_1, r_2;n)q^n.\eeq 
Our main result is the following. See Section \ref{sec:gen_functions} for a more detailed discussion of the cohomological invariants for which we state exact formulas and asymptotics below.


\begin{theorem}\label{thm:main_result}
Let $S$ be a smooth projective surface.
Then we have the following exact formulas:
\ben 
\item If $0 \leq \chi(S) < 24n$, then we have \label{eq:Euler_char_exact_formula}
\begin{align*} 
    \chi(\mathrm{Hilb}^n(S))=&2\pi
    \sum_{j<{\frac{\chi(S)}{24}}}
    \sum_{{k=1}}^\infty
    k^{\chi(S)/2}A_{k}(-\chi(S),0,j;n)
    \chi(\mathrm{Hilb}^j(S))
    L^*(0,j,k;n).
\end{align*}
\item If $\sigma(S) \leq \chi(S)< 24n$, then we have
\begin{align*}
     \sigma(\mathrm{Hilb}^n(S))= &
      \hspace{5mm} 2\pi
    \sum_{j<{\frac{\chi(S)}{24}}}
    \sum_{\substack{k=2\\k\mathrm{\ even}}}^\infty
    \frac{A_{k}(\sigma(S),\Lambda(S),j;n)
    \sigma(\mathrm{Hilb}^j(S))
    }{k^{\Lambda'(S)/2}}L^*(0,j,k;n)\\ \nonumber
    &+
    2\pi\sum_{j<\frac{3\sigma(S) -\chi(S)}{48}}
     \sum_{\substack{k=1\\k \ \mathrm{ odd}}}^\infty \frac{(-1)^nB_{k}(\sigma(S),\Lambda(S),j;n)a(\Lambda(S),\sigma(S);j)}{ 2^{\Lambda(S)/2}k^{\Lambda'(S)/2}}L^*(1,j,k;n)
     .
\end{align*}
where $\Lambda(S):=-(\sigma(S) + \chi(S))/2$ and $\Lambda'(S):=(\sigma(S) - \chi(S))/2$.
 \een
Here $a(\alpha,\beta,j)$, $A_{k}(\alpha,\beta,j;n)$, $B_{k}(\alpha,\beta,j;n)$,   $L^*(0,j,k;n)$, and $L^*(1,j,k;n)$ are defined in Section \ref{sec:Circle_method_outline}, (\ref{eq:Kloosterman_def1}), (\ref{eq:Kloosterman_def2}),   (\ref{eq:L0Bessel}), and (\ref{eq:L1Bessel}), respectively. 
 \end{theorem}

{ \noindent \tbf{Remark:}
In Theorem \ref{thm:main_result}, $A_k(\alpha,\beta,j;n)$ and $B_{k}(\alpha,\beta,j;n)$ are known as Kloosterman sums, and $L^*(0,j,k;n)$ and $L^*(1,j,k;n)$ are, up to simple multiplicative factors, modified Bessel functions of the first kind. \medskip
}

Theorem \ref{thm:main_result} offers exact formulas as convergent infinite series. These formulas  imply the following asymptotics.

\begin{corollary}\label{cor:asymptotics}
Let $S$ be a smooth projective surface. 
Then the following are true:
\begin{enumerate}
    \item  Suppose $\chi(S)  \geq \sigma(S)$.
    \bena 
     \item
    If $\sigma(S) < 0,$ then as $n \to \infty,$ we have 
    \bal \sigma(\mathrm{Hilb^n}(S)) 
\sim &
(-1)^n 2^{\frac{7\sigma(S)-3\chi(S)-14}{8}}3^{\frac{\sigma(S)-\chi(S)-2}{8}}
(\chi(S)-3\sigma(S))^{\frac{\chi(S)-\sigma(S)+2}{8}}
n^{\frac{\sigma(S)-\chi(S)-6}{8}} \\
& 
\cdot \exp\left(\pi\sqrt{\frac{\chi(S)-3\sigma(S)}{6}n}\right). \eal 
    \item 
    If $\sigma(S) >0$, then as $n \to \infty,$ we have 
    
    $$\sigma(\mathrm{Hilb^n}(S)) \sim 
2^{\frac{3\sigma(S) -3\chi(S)-14}{8}}
3^{\frac{\sigma(S)-\chi(S)-2}{8}}
\chi(S)^{\frac{ \chi(S) -\sigma(S)+2}{8}}
n^{\frac{\sigma(S)-\chi(S)-6}{8}}
\exp\left(\pi\sqrt{\frac{\chi(S)}{6}n}\right).$$

\item
    If $\sigma(S)=0$ and $\chi(S) \neq 0$, then as $n\to\infty$, we have 
    $$\sigma(\mathrm{Hilb^{2n}}(S)) \sim  2^{\frac{-\chi(S) -3}{2}}
3^{\frac{-\chi(S)-2}{8}}
\chi(S)^{\frac{2+\chi(S)}{8}}
n^{\frac{-\chi(S)-6}{8}}
\exp\left(\pi\sqrt{\frac{\chi(S)}{3}n}\right). $$ Moreover, for all $n$ we have
$$\sigma(\mathrm{Hilb^{2n+1}}(S)) = 0.$$
\item
If $\sigma(S) =\chi(S) =0$, we have for all $n$,
$$\sigma(\mathrm{Hilb^{n}}(S)) = 0.$$
    
    \een
    
    \item \bena \item If $\chi(S) >  0,$ then, as $n \to \infty,$ we have $$\chi(\mathrm{Hilb}^n(S)) \sim 2^{\frac{-3\chi(S)-5}{4}}3^{\frac{-\chi(S)-1}{4}}\chi(S)^{\frac{\chi(S)+1}{4}}n^{\frac{-
\chi(S)-3}{4}}\exp\left(\pi\sqrt{\frac{2\chi(S)}{3}n}\right).$$
\item If $\chi(S) =0$, then for all $n$ we have
$$\chi(\mathrm{Hilb}^n(S)) =0.$$
\een
\end{enumerate}
\end{corollary}

The asymptotics in Corollary 1.2 imply the following asymptotic  properties of $b^*_S(r;n)$ and $c_S^*(r_1,r_2;n)$.

\begin{corollary}\label{cor:equidistribution}
Let $S$ be a smooth projective surface. 
Then the following are true:
\begin{enumerate}
    \item Suppose that  $\chi(S) \geq \sigma(S)$. If  $\chi(S)+\sigma(S)>0$, then as $n \to \infty $ we have
    $$b^*_S(0;n) \sim b^*_S(1;n). $$
    If $\chi(S)+\sigma(S) =0$, then as $n \to \infty $ we have
    $$ b^*_S(1;n) =0.
    $$ 
    If $\chi(S)+\sigma(S) <0$, then as $n \to \infty $ we have
    $$b^*_S(0;n) \sim - b^*_S(1;n). $$
    \item  Suppose that $\chi(S) \geq \sigma(S)$. If $h^{1,0}=0$, then as $n \to \infty $ we have
    $$c_S^*(0,0;n) \sim c_S^*(1,1;n) \ \ \text{and} \ \  c_S^*(0,1;n)=c_S^*(0,1;n) =0.$$
     If  $h^{1,0}>0$, then as $n \to \infty $ we have
    $$c^*_{S}(0,0;n) \sim c^*_{S}(1,1;n) \sim -c^*_{S}(0,1;n) = -c^*_{S}(1,0;n).$$
\end{enumerate}
\end{corollary}
Note that when $S$ is a K3 surface, Corollary \ref{cor:equidistribution} (1) recovers the equidistribution of the $b_S^*(r; n)$ that follows from the work of  Manschot and Zapata Rolon in \cite{Manschot}.
\medskip

\noindent \tbf{Remark:} The Enriques-Kodaira Classification Theorem \cite[p.~244]{CCS} mostly describes the possible Hodge structures of smooth complex surfaces. This theorem lists minimal models for many birational equivalence classes of surfaces, which determine $h^{0,0}$, $h^{1,0}$, $h^{2,0}$,  and $\min\{h^{1,1}\}$ for that class. By the blowup construction, every smooth complex projective surface $S$ is birationally equivalent to a smooth complex projective surface $S'$ with $h^{1,1}(S')=h^{1,1}(S)+1$ and $h^{s,t}(S') = h^{s,t}(S)$ for $(s,t) \neq (1,1)$. A minimal model of a birational equivalence class is a surface that is not a blowup of any other smooth surface. Since $\chi(S)$ and $\chi(S)-\sigma(S)$ increase linearly with $h^{1,1}(S)$, Theorem \ref{thm:main_result}, Corollary \ref{cor:asymptotics}, and Corollary \ref{cor:equidistribution} apply to all surfaces in each birational equivalence class except those whose Hodge structures fall in a certain finite set. In particular, if the minimal model satisfies the hypotheses of these statements, then  all surfaces in that class satisfy them. 
Excluding surfaces of general type, the only classes of projective surfaces whose minimal models do not satisfy these hypotheses are
ruled surfaces of genus $g \geq 2$ 
(see \cite{CCS} and
\cite{Schutt}).

\skipaline

In Section \ref{sec:gen_functions}, we will define all of the terms above and describe their topological and geometric significance. We will also state important properties of the generating functions of these sequences, which will be crucial in our deduction of the above formulas. In Section \ref{sec:Circle_method_outline} we will outline { our use of the circle method to prove} Theorem \ref{thm:exactFormulasFor_H}, a general result from which Theorem \ref{thm:main_result} and Corollaries \ref{cor:asymptotics} and \ref{cor:equidistribution} are derived. { Sections \ref{sec:even_case} and} \ref{sec:odd_case} include arguments necessary for this proof. In Section \ref{sec:proofs}, these arguments are assembled, and Theorem \ref{thm:exactFormulasFor_H}, Theorem \ref{thm:main_result}, and Corollaries \ref{cor:asymptotics} and \ref{cor:equidistribution} are proven. In Section \ref{sec:numerics} we illustrate our results with numerics. 

\section{{ Preliminaries}} \label{sec:gen_functions}
{ In this section, we present G\"ottsche's result and specialize it in terms of weakly holomorphic modular forms. 
We then provide bounds on certain exponential sums known as Kloosterman sums and bounds on $I$-Bessel functions. Both Kloosterman sums and Bessel functions will appear in our application of the circle method in Sections \ref{sec:even_case} and \ref{sec:odd_case}.}

A compact complex manifold $M$ has cohomological invariants called the Hodge numbers $h^{s,t}:=h^{s,t}(M),$ which are defined as the complex dimensions of the $(s,t)$-Dolbeault cohomology space $H^{s,t}(M)$ (see \cite{Wells}). 
When the context is clear, we will not explicitly indicate dependence of the Hodge numbers $h^{s,t}$ on the manifold $M$.
If $M$ is a K\"ahler manifold, the Hodge numbers are related to the Betti numbers $b_n(M)$ by the formula
\beq b_n(M)= \sum_{s+t=n}h^{s,t}(M) \label{eq:Betti_Hodge} \eeq
(see \cite[p.~198]{Wells}). Moreover, for any manifold $M$ one can construct a manifold $\text{Hilb}^n(M)$ which can be thought of as a smoothed version of the $n^{th}$ symmetric product of $M$ (see \cite{Hart}). For any smooth projective complex surface $S,$ G\"ottsche's formula allows one to compute the Hodge numbers of  $\text{Hilb}^n(S)$ for all $n$ from the Hodge numbers of $S$.

In order to state G\"ottsche's result, we first define the Hodge polynomial, which serves as a generating function for the Hodge numbers of $M$:
\beq \begin{aligned}
 \label{eq:chihodge} \chi_{\mathrm{Hodge}}(M)(x,y)
:= & x^{-d/2}y^{-d/2}\sum_{s,t}h^{s,t}(M)(-x)^s(-y)^t, 
\end{aligned} \eeq 
where $d$ is the complex dimension of $M$. We will generally supress the $(x,y)$ for notational convenience.


On page 37 of \cite{GotBook}, G\"ottsche proved the remarkable fact that one can assemble the Hodge polynomial for $\mathrm{Hilb}^n(S)$ using the Hodge numbers $h^{s,t}(S)$.

\bthm[G\"ottsche]\label{thm:gottsche} If $S$ is a smooth projective complex surface, then we have that
$$Z_S(x,y;\tau):=\sum_{n\geq0}\chi_{\mathrm{Hodge}}(\mathrm{Hilb}^n(S))q^n =\prod_{n =1}^\infty \frac{\prod_{s +t \mathrm{\ odd}} (1- x^{s-1}y^{t-1}q^n)^{h^{s,t}}}{\prod_{s +t \mathrm{\ even}} (1- x^{s-1}y^{t-1}q^n)^{h^{s,t}}}. $$
\ethm
\noindent \tbf{Remark:}
Note that for each $n$, $|s|,|t|\leq n$ in the definition of $\chi_{\mathrm{Hodge}}(\mathrm{Hilb}^n(S))$.\medskip

\noindent By \cite[p.~43]{Ball}, every such smooth projective complex surface $S$ is K\"ahler, so by Serre duality and Hodge symmetry, the Hodge numbers satisfy the following relations:
$$ h^{0,0}=h^{2,2} \qquad \qquad  h^{1,0}=h^{0,1}= h^{1,2}=h^{2,1} \qquad  \qquad h^{2,0}=h^{0,2}. $$
By the additivity of the Hodge numbers, we need only consider the case where $M$ is connected, i.e. $h^{0,0}=1$. In this case, we obtain:
\beq \label{eq2} Z_S(x,y;\tau)=\prod_{n =1}^\infty \frac{\big((1- x^{-1}q^n)(1- xq^n)(1- y^{-1}q^n)(1-yq^n)\big)^{h^{1,0}}}{(1- x^{-1}y^{-1}q^n)(1- xyq^n)\big((1- x^{-1}yq^n)(1- xy^{-1}q^n)\big)^{h^{2,0}}(1-q^n)^{h^{1,1}}}. \eeq

This polynomial gives rise to other topological invariants upon substituting $\pm 1$ for $x$ and $y.$ For example, the Hirzebruch $\chi_y$-genus of $M$ is the polynomial
\beq \label{eq:chiy} \chi_y(M):=\sum_{s,t} (-1)^t h^{s,t}(M) y^s. \eeq
We can express this in terms of the Hodge polynomial:
$$ \chi_{\text{Hodge}}(\text{Hilb}^n(S))(y,1)=  y^{-n} \sum_{s,t} (-1)^t h^{s,t}(\text{Hilb}^n(S))(-y)^s =  \chi_{-y}(\text{Hilb}^n(S))y^{-n}. $$
In terms of the Betti numbers, the Euler characteristic $\chi(M)$ is defined as
\beq \label{euchar} \chi(M):= \sum_n (-1)^n b_n(M). \eeq
Setting $x = 1$ and $y = 1$ in the Hodge polynomial, we see in reference to (\ref{eq:Betti_Hodge}) and (\ref{euchar}) that  
\beq \label{eq:euch} \chi_{\text{Hodge}}(\text{Hilb}^n(M))(1,1) =  \chi_{-1}(\text{Hilb}^n(M))=  \chi(\text{Hilb}^n(M)). \eeq
On the other hand, setting $x = -1$ and $y = 1$ in the Hodge polynomial gives
\beq \label{eq:sig} \chi_{\text{Hodge}}(\text{Hilb}^n(M))(-1,1) =  (-1)^n \chi_{1}(\text{Hilb}^n(M))=(-1)^n  \sigma(\text{Hilb}^n(S)),\eeq where the signature $\sigma(M)$ of a $d$-dimensional complex manifold $M$ is the signature of the intersection pairing on $H^d(M)$ (see \cite{May}).
In terms of Hodge numbers of K\"ahler surfaces $S$, the signature is given by
$$\sigma(S)= 2h^{2,0}+2-h^{1,1}.$$

The discussion above indicates the importance of G\"ottsche's infinite product formulas (Theorem \ref{thm:gottsche}) as a vehicle for studying  invariants of complex projective surfaces. This is further illuminated by the following specializations of Theorem \ref{thm:gottsche} in terms of these  invariants.

\begin{lemma}\label{lem:gottsche_prod} If $S$ is a smooth projective complex surface, then the following are true:
\begin{align} \sum_{n =0}^\infty \chi_{-y}(\mathrm{Hilb}^n(S))y^{-n}q^n= & \prod_{n =1}^\infty \frac{((1- y^{-1}q^n)(1-yq^n))^{h^{1,0}-h^{2,0}-1}}{(1-q^n)^{h^{1,1}-2h^{1,0}}}, \label{sp1} \\
\sum_{n=0}^\infty \chi(\mathrm{Hilb}^n(S))q^n= & \prod_{n=1}^\infty (1-q^n)^{-\chi(S)}, \label{sp2} \\
\sum_{n=0}^\infty  (-1)^n \sigma(\mathrm{Hilb}^n(S))q^n= & \prod_{n=1}^\infty \frac{(1-q^n)^{\sigma(S)}}{(1-q^{2n})^{(\sigma(S) + \chi(S))/2}}. \label{sp3}
\end{align} 
\end{lemma}

\noindent \tbf{Remark:} We note that (\ref{sp2}) and (\ref{sp3}) are alternate expressions for $Z(1,1;\tau)$ and $Z(1,-1;\tau)=Z(-1,1;\tau)$, respectively. We can apply the same process to $Z(-1,-1;\tau)$ to obtain
\beq \label{eq:minus_expression} Z(-1,-1;\tau)= \prod_{n=1}^\infty \frac{(1-q^n)^{4h^{1,0}}}{(1-q^{2n})^{\chi(S)+8h^{1,0}}}.
\eeq 

We now show that these functions can be assembled in linear combinations to give alternate formulas for $B_S(r;\tau)$ and $C_S(r_1,r_2;\tau)$.



\blem\label{lem:roots_of_unity} 

Let $S$ be a smooth projective complex surface. We have
\beq \label{eq:B_sum} B_S(r, \ell;n) = \frac{1}{2}\left( Z_S(1,1;\tau) + (-1)^r Z_S(1,-1;\tau)\right) \eeq
and
\beq  
 \label{eq:now} C_S(r_1, r_2;\tau) = \frac{1}{4} \sum_{\substack{ j_1 \Mod{2} \\ j_2\Mod{2}}}  (-1)^{j_2r_2} (-1)^{j_1r_1} Z_S((-1)^{j_2}, (-1)^{j_1};\tau), \eeq
 where $C_S(r_1, r_2;\tau)$ is defined by (\ref{eq:bigC}). \elem

\bpf We prove only (\ref{eq:now}), since (\ref{eq:B_sum}) { follows from even simpler manipulations.} 
We have 
\bal  
\sum_{n \geq 0} c^*_S(r_1, r_2;n) q^n 
= & \frac{1}{4} \sum_{s, t ,n } c_S(s,t;n)\sum_{\substack{ j_1 \Mod{2} \\ j_2 \Mod{2}}} (-1)^{j_2(s+r_2)} (-1)^{j_1(t+r_1)}q^n  \\
= & \frac{1}{4} \sum_{\substack{ j_1 \Mod{2} \\ j_2\Mod{2}}}(-1)^{j_2r_2}  (-1)^{j_1r_1} \sum_{s, t , n} c_S(s, t; n) (-1)^{j_2s} (-1)^{j_1t}q^n  \\
= & \frac{1}{4} \sum_{\substack{ j_1\Mod{2} \\ j_2\Mod{2}}} (-1)^{j_2r_2} (-1)^{j_1r_1} Z_S((-1)^{j_2}, (-1)^{j_1};\tau).
\eal \epf

\subsection{Modularity Properties of Specializations of G\"ottsche's Identity \label{sec:trans_law}}

The infinite products 
(\ref{sp2}), (\ref{sp3}), and (\ref{eq:minus_expression}) can  be written in terms of the Dedekind eta function, the weight 1/2 modular form { on $SL_2(\bZ)$ (with multiplier)} that is defined by 
\begin{equation}\label{eq:eta}
    \eta(\tau) := q^{1/24} \prod_{n=1}^\infty (1-q^n).
\end{equation} This is done explicitly in the following lemma. 

 
 
\begin{lemma}\label{lem:z_to_h}
Let $S$ be a smooth projective complex surface.  We have
 that
 \begin{align}
  \label{sp5} Z_S(1,1;\tau) = &  q^{{\chi(S)}/{24}} \frac{1}{{\eta(\tau)}^{\chi(S)}}, \\
  \label{sp6}  Z_S(-1,1;\tau) = Z_S(1,-1;\tau)  = &  q^{{\chi(S)}/{24}}\frac{\eta(\tau)^{\sigma(S)}}{\eta(2 \tau)^{(\sigma(S) +\chi(S))/{2}}}, \\
  \label{sp7}  Z_S(-1,-1;\tau) = & q^{\chi(S)/24}\frac{\eta(2 \tau)^{4h^{1,0}}}{\eta(\tau)^{\chi(S)+8h^{1,0}}}.
 \end{align}

\end{lemma}

{ Having expressed $Z_S( \pm 1, \pm1; \tau)$ as an eta quotient (up to a fractional power of $q$)}, we define
\begin{equation}\label{eq:hab_def}
    H_{\alpha, \beta}({ q}) := q^{\chi(S)/24} \eta(\tau)^\alpha \eta(2 \tau)^\beta =: Z_S( x,y; \tau),
\end{equation}
where $x,y = \pm 1$ and $\alpha$ and $\beta$ are determined by $x, $ $y$, and  the Hodge numbers of $S$. {  Thus, (up to a fractional power of $q$)}, $H_{\alpha, \beta}$ is a modular form on $\Gamma_0(2)$. { In our application of the circle method} to prove our main results, we  use the modularity of $H_{\alpha,\beta}$ to determine the behavior of this function when $\tau$ is near the rational number $h/k$. 
The modular transformation equations for $H_{\alpha, \beta}$ are given in the following lemma.
{
The $q$-expansion of $H_{\alpha,\beta}$ gives the behavior of $H_{\alpha,\beta}$ for $\tau$ near  $i \infty$, and this lemma relates the behavior of $H_{\alpha,\beta}$ near  $i\infty$ to its behavior near  $h/k.$
We present two transformation equations, corresponding to the two different cusps of $\Gamma_0(2)$.
These cusps 
partition the rationals $h/k$ into two parts, based on the parity of $k$. The two cusps will  make distinct contributions  to our final formula.  
}

\blem\label{lem:trans_law} Let $k$ be a positive integer, let $h$ satisfy the condition $(h,k) =1$, and let $h'$ satisfy $hh' \equiv -1 \mod k$. 
Further suppose $\mathrm{Re}(z) >0$.
Define
$$\omega_{\alpha, \beta}(h,k) :=\exp\left(- \pi i( \alpha s(h,k)+\beta s(2h,k))\right),$$
where $s(h,k)$ is the Dedekind sum
\bal s(h,k) := &\sum_{r=1}^{k-1} \bigg(\bigg(\frac{r}{k}\bigg)\bigg)\bigg(\bigg( \frac{hr}{k}\bigg)\bigg), \eal
and $((x))$ is the sawtooth function
\bal
((x)) := & \bcas 0 & x \in \bZ \\
x-\lfloor x\rfloor  - \frac{1}{2} & x \notin \bZ. \ecas \eal
\begin{enumerate}
    \item If $2 \mid k$, we have
\begin{align*}  H_{\alpha,\beta} &\left( \exp \left( \dfrac{2 \pi i h}{k} - \dfrac{2 \pi z}{k} \right) \right) \\
    &=    
     z^{-\frac{\alpha + \beta}{2}}\omega_{\alpha, \beta}(h,k)
     \cdot \exp \left( (\alpha + 2 \beta)\dfrac{\pi}{12 k} \left( z -\dfrac{1}{z} \right) \right) 
     \cdot H_{\alpha, \beta} \left( \exp \left( \dfrac{2 \pi i h'}{k} - \dfrac{2 \pi}{z k} \right) \right).
\end{align*}
\item If $2 \nmid k$, we have
 \begin{align*} 
      H_{\alpha, \beta} & \left(\exp\left( \dfrac{2 \pi i h}{k} - \dfrac{2 \pi z}{k} \right)\right)\\  \hspace{5mm} \qquad & =  \frac{z^{-\frac{\alpha + \beta}{2}}\omega_{\alpha, \beta}(h,k)}{2^{\beta/2}}
     \cdot \exp \left( \frac{\pi}{24k}\left( {2(\alpha + 2 \beta)}z -\frac{2\alpha + \beta}{z}\right) \right)  \cdot H_{\beta,\alpha}\left(\exp\left(\frac{\pi ih'}{k}  - \frac{\pi}{zk} \right)\right).
 \end{align*}
\end{enumerate}
 
\elem
\bpf 
From \cite[p.~96]{Apostol}, we know 
{\begin{align*}
    H_{-1,0}&\left(\exp\left(\dfrac{2 \pi i h}{k} - \dfrac{2 \pi z}{k}\right)\right)\\ &= z^{\frac{1}{2}} \exp(\pi i s(h,k)) \exp\left(\dfrac{\pi}{12k} \left(\dfrac{1}{z} - z\right)\right) H_{-1,0}\left(\exp\left(\dfrac{2 \pi i h'}{k} - \dfrac{2 \pi}{k z}\right)\right),
\end{align*}
where $(h,k) =1$, ${\tau= (iz + h)/k}$, $hh' \equiv -1 \mod k$, and $s(h,k)$ is the Dedekind sum defined above. Similarly, \cite{Hagis3} shows that
if $2 \mid k$, then
\begin{align*}
    H_{-1,1}&\left(\exp\left( \dfrac{2 \pi i h}{k} - \dfrac{2 \pi z}{k} \right)\right)\\ &= \exp ( \pi i (\sigma(h,k) ) \exp \left( \dfrac{\pi}{12 k} \left(z - \dfrac{1}{z} \right) \right) H_{-1,1} \left( \exp \left( \dfrac{2 \pi i h'}{k} - \dfrac{2 \pi}{z k} \right) \right), 
\end{align*}
where $\sigma(h,k) := s(h,k)-s(2h,k)$, and if $2 \nmid k$
\begin{align*}
    H_{-1,1}&\left(\exp\left( \dfrac{2 \pi i h}{k} - \dfrac{2 \pi z}{k} \right)\right)\\ &= 2^{-\frac{1}{2}}\exp(\pi i \sigma(h,k))\exp\left(\frac{\pi}{24k}\left(\frac{1}{z} + 2z\right)\right)H_{-1,1}\left(\exp\left(\frac{ \pi ih'}{k}  - \frac{\pi}{zk} \right)\right)\inv
\end{align*}}
\noindent To prove the lemma, { we}  { note} that $H_{\alpha, \beta}({ q}) = H_{-1,0}({ q}) ^ {-(\alpha + \beta)} \cdot H_{-1,1}({ q}) ^ \beta$ and $H_{\alpha,0} ({ q^2}) = H_{0,\alpha}({ q})$. 
\epf


{ \noindent \tbf{Remark:}
When we apply the circle method to find the coefficients of the $q$-expansion of $H_{\alpha, \beta}$ in Sections \ref{sec:Circle_method_outline}, \ref{sec:even_case}, and \ref{sec:odd_case}, we consider $\tau \in \bH$ on the horizontal line segment $$L = \{\tau \ | \ \tau = u + N^{-2} i, \  0 \leq u \leq 1 \}.$$ As $N \to \infty$, 
$L$ approaches the positive real axis directly from above. As discussed in Section \ref{sec:Circle_method_outline}, for the purpose of estimates we partition $L$ into line segments that have midpoints at $h/k + i N^{-2}$ for cusp representatives $0 \leq h / k  \leq 1$ and $k \leq N.$
When we reparameterize these line segments in terms of the variable $z$ introduced in the statement of Lemma \ref{lem:trans_law}, they become small vertical line segments just to the right of the origin in the $z$-plane. Lemma \ref{lem:trans_law} allows us to estimate the value of $H_{\alpha, \beta}$ for $z$ on these small vertical line segments. 
}

{
\subsection{Bounds on certain Kloosterman sums\label{sec:Kloosterman}}

In order to bound the error terms that will result from estimating the contour integral \eqref{eq:countour_integral}, we need bounds on Kloosterman sums of the form
\begin{align} \label{eq:Kloosterman_def1} A_k(\alpha,\beta,j;n):= & \sum_{\substack{0\leq h<k\\(h,k)=1}}\omega_{\alpha, \beta}(h,k)\exp\left(-\frac{2\pi inh}{k}+\frac{2\pi ih'j}{k}\right) \\
=: & \sum_{\substack{0\leq h<k\\(h,k)=1}} A_{h,k}(\alpha,\beta,j;n)
\end{align}
and 
\begin{align} \label{eq:Kloosterman_def2} B_k(\alpha,\beta,j;n):= & \sum_{\substack{0\leq h<k\\(h,k)=1}}\omega_{\alpha, \beta}(h,k)\exp\left(-\frac{2\pi inh}{k}+\frac{\pi ih'j}{k}\right) \\
=: & \sum_{\substack{0\leq h<k\\(h,k)=1}} B_{h,k}(\alpha,\beta,j;n).
\end{align}
These sums admit the trivial estimate $O(k)$. However, in the case $\alpha + \beta = 0$, we will need a sharper estimate. In addition, we will occasionally restrict the values of $k_i$ to an interval $N-k < k_i \leq \sigma < N$, which will in turn restrict $h'$ to one or two intervals modulo $k$. For bounding purposes, it suffices to consider sums of the form 
$$ \sideset{}{'}\sum_{\substack{0\leq h<k\\(h,k)=1}} A_{h,k}(\alpha,\beta,j;n) \qand \sideset{}{'}\sum_{\substack{0\leq h<k\\(h,k)=1}} B_{h,k}(\alpha,\beta,j;n)$$
where the $'$ indicates that $h'$ is restricted to an interval $0 \leq \sigma_1 \leq h' < \sigma_2 \leq k$. 
Thus we will need the following lemma.
\blem For $\alpha$ and $\beta$ fixed, $\alpha + \beta =0$, the sums 
$$ \sideset{}{'}\sum_{\substack{0\leq h<k\\(h,k)=1}} A_{h,k}(\alpha,\beta,j;n) \qand \sideset{}{'}\sum_{\substack{0\leq h<k\\(h,k)=1}} B_{h,k}(\alpha,\beta,j;n)$$
are subject to the estimate 
$$O(n^{1/3}k^{2/3 + \ep}) $$
uniformly in $\sigma_1$, $\sigma_2$, and $j$. \label{lem:Kloosterman}
\elem
\bpf
The proof is a simple adaptation of those of Theorems 2 and 3 of \cite{Hagis3}, which follow the proof of Theorem 2 in \cite{Lehner}. Equations (3.5) and (4.11) in \cite{Hagis3} together state that if $2 \mid k$, we have
$$\exp({ \pi i (s(h,k)-s(2h,k)})) = \exp\left(2 \pi i \left( \frac{4\phi(uh + vh')}{(k/2,2)Gk}+r(h,k)\right)\right) $$ where $r(h,k)$ is a rational number that depends on $k$ and $h $ mod 4, $\phi$, $G$, and $v$ are integers that depend only on $k$, and $u$ is a polynomial in $k$. (5.7) states that if $2 \nmid k$, then
$$\exp({ \pi i (s(h,k)-s(2h,k)})) = \exp\left(2 \pi i \left( \frac{\Phi(uh + vh')}{gk}+r(h,k)\right)\right) $$
with $r$, $v$ and $u$ as before, and $\Phi$ and $g$ integers dependent on $k$. The proof of Lemma \ref{lem:Kloosterman} proceeds as in \cite{Hagis3}, with $u$, $v$, and $r$ replaced with $\beta u$, $\beta v$ and $\beta r$, respectively. \epf
}

{  \subsection{Bounds on $I$-Bessel functions}
As noted in the introduction, an important ingredient in our exact formulas is the modified Bessel function of the first kind, which is also known as the $I$-Bessel function. Here we recall a number of facts from \cite{NIST} about $I$-Bessel functions that will be useful throughout our arguments. We define the $I$-Bessel function of order $v$ as
\begin{equation*}
    I_v(z) := \bigg( \dfrac{z}{2} \bigg)^v \sum_{k = 0}^\infty \dfrac{\big( \frac{1}{4} z^2 \big)^k}{k! \Gamma(v + k + 1)}.
\end{equation*}
From this definition it follows that for $0 < x <  1,$ we have
\begin{equation}\label{eq:geom_bes_bound}
    |I_v(x)| < \frac{4}{3} \cdot \Big( \frac{x}{2} \Big)^v.
\end{equation}
We also recall an integral representation of the $I$-Bessel function which we will use in Sections \ref{sec:even_case} and \ref{sec:odd_case}: for $v \geq 0,$ we have
\begin{equation}\label{eq:bes_int}
    I_v(z)=\frac{\left(\frac{1}{2}z\right)^{-v}}{2\pi i}\int_Rt^{v-1}\exp\left[t+\frac{z^2}{4t}\right]dt,
\end{equation}
where $R$ is any simple closed contour surrounding the origin. 
}

\section{Outline of the circle method \label{sec:Circle_method_outline}}

In this section, we outline our use of the circle method to find exact formulas for coefficients $a(\alpha, \beta;n)$ of $H_{\alpha, \beta}$ 
where $\alpha+\beta\leq 0$.
This corresponds to those cases where $H_{\alpha,\beta}$ is essentially a modular form of non-positive weight.
{ These cases are special because non-constant
holomorphic modular forms do not exist in non-positive weight.} 
We make use of the method of Rademacher, introduced in \cite{Rademacher_partitions}, which improved upon the earlier work of G. H. Hardy and Srinivasa Ramanujan on the partition function. We adapt the implementation of this method along the lines of earlier works by Hagis (see \cite{Hagis3}, \cite{Hagis4}, and \cite{Hagis5}).

The basic ingredient of the circle method is Cauchy's integral formula, which we use to obtain
\begin{equation}\label{eq:countour_integral}
a(\alpha,\beta;n)=\frac{1}{2\pi i}\int_C \frac{H_{\alpha, \beta}(q)}{q^{n+1}}dq.
\end{equation}
{ We let $C$ be the circle of radius $e^{-2\pi N^{-2}}$ centered at the origin, where $N$ is an upper bound on the denominator of the cusp representatives under consideration. 
We emphasize that we choose $C$ to be a circle with this particular radius; there are many paths that could be chosen in applying Cauchy's integral formula, but the one we choose here is a particular choice following the implementation of the circle method in \cite{Hagis3}, \cite{Hagis4}, and \cite{Hagis5}.}
The exponent $n$ will be fixed, and $N$ will later be allowed to approach infinity. As $N \to \infty$, the main contribution to this integral comes from a dense set of poles on the boundary of the unit circle, where each pole is located at a root of unity $e^{2 \pi i h/k}$. 
{ The locations of these poles follow from the definition of $H_{\alpha, \beta}$ in \eqref{eq:hab_def} and the location of the zeros of $\eta(\tau)$ afforded by the product representation \eqref{eq:eta}.}

The generating function  $H_{\alpha, \beta}$ is,  up to a fractional power of $q$, a  modular form with non-positive weight on $\Gamma_0(2)$ in the cases we consider. As shown in Lemma \ref{lem:trans_law}, we have that $H_{\alpha, \beta}$ has a pole of order $-(\alpha + 2\beta)/24$ at the cusp of $\Gamma_0(2)$ whose representatives have all even denominators and a pole of order $-(2\alpha + \beta)/48$ at the cusp whose representatives have all odd denominators.

\begin{figure}[h]
\centering
\includegraphics[scale=.5]{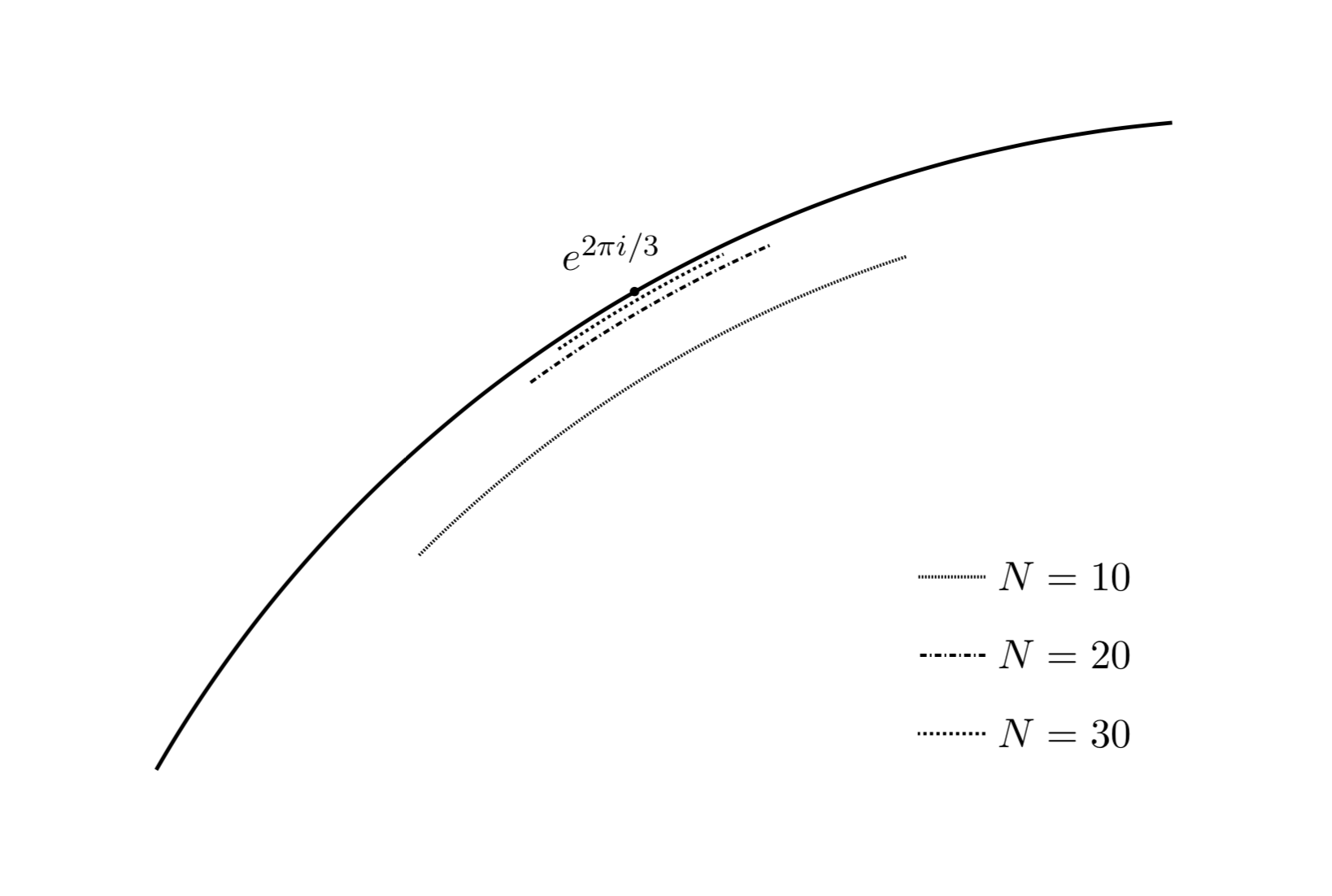}
\caption{Farey arcs $\xi_{1,3}$}
\label{fig:Farey}
\end{figure}

Since we have descriptions of $H_{\alpha,\beta}$ near all roots of unity, we divide $C$ into \emph{Farey arcs} $\xi_{h,k}$ for $(h,k)=1$, which decrease in length and  approach the point $e^{2 \pi i h/k}$ as $N$ increases. 
Figure \ref{fig:Farey} shows the Farey arc $\xi_{1,3}$ for $N=10$, $20$, and $30$.
By estimating $H_{\alpha,\beta}$ along these Farey arcs, we obtain a convergent series for the coefficients $a(\alpha; \beta, n)$ as $N \to \infty.$ 
These series are exact formulas which yield Theorem \ref{thm:exactFormulasFor_H}.

We will need explicit descriptions of these Farey arcs.
First, we can rewrite
\[\frac{1}{2\pi i}\int_C \frac{H_{\alpha, \beta}(q)}{q^{n+1}}dq=\sum_{\substack{0\leq h<k\leq N\\(h,k)=1}}\frac{1}{2\pi i}\int_{\xi_{h,k}} \frac{H_{\alpha, \beta}(q)}{q^{n+1}}dq.\]
In the Farey series of order $N$, we consider the fraction $h/k$ and its two neighbors,
\beq \frac{h_1}{k_1}<\frac{h}{k}<\frac{h_2}{k_2}, \label{eq:neighbors} \eeq
discussed in \cite{Rademacher}. Each of $h_i$ and $k_i$ depend on $h$, $k$ and $N$ as described in \cite{Rademacher}. Note that $N-k <  k_i \leq  N$.
{ On the arc $\xi_{h,k}$, we can introduce the variable $\theta$ via the  transformation
\[q=\exp\left[-2\pi N^{-2}+2\pi i\left(\frac{h}{k}+\theta\right)\right]=\exp\left[\frac{2\pi ih}{k}-\frac{2\pi z}{k}\right],\]}
where 
\[-\vartheta_{h,k}':=-\frac{1}{k(k_1+k)} \leq \theta \leq  \frac{1}{k(k_2+k)}=: \vartheta_{h,k}''.\]
Setting $z=k(N^{-2}-i\theta)$ on each arc $\xi_{h,k}$, we obtain
\begin{equation*}
  a(\alpha,\beta;n)
=\sum_{\substack{0\leq h<k\leq N\\(h,k)=1}}e^{-\frac{2\pi inh}{k}}\int_{-\vartheta_{h,k}'}^{\vartheta_{h,k}''}H_{\alpha,\beta}\left(e^{\frac{2\pi ih}{k}-\frac{2\pi z}{k}}\right)e^{\frac{2\pi nz}{k}}d\theta.
\end{equation*}
Guided by Lemma \ref{lem:trans_law}, we break up the sum into even and odd $k$, writing
\begin{equation}
    a(\alpha,\beta;n)=S(0,N;n)+S(1,N;n),\label{eq:a_parity_strat_S's}
\end{equation}
where
$$S(r,N;n) := \sum_{\substack{0\leq h<k\leq N\\(h,k)=1\\ k \equiv r \Mod{2}}}e^{-\frac{2\pi inh}{k}}\int_{-\vartheta_{h,k}'}^{\vartheta_{h,k}''}H_{\alpha,\beta}\left(e^{\frac{2\pi ih}{k}-\frac{2\pi z}{k}}\right)e^{\frac{2\pi nz}{k}}d\theta.
$$
We henceforth omit dependence on $\alpha$ and $\beta$ from the names of most variables.
Section \ref{sec:Kloosterman} gives estimates on sums of roots of unity known as Kloosterman sums that are needed to bound the error terms.
Section \ref{sec:even_case} extracts behavior from $S(0,N;n)$ that will contribute to the exact formula and bounds the error terms. 
Throughout the following three sections, $n$, $\alpha$, and $\beta$ are fixed.

\section{The even case}\label{sec:even_case}

\subsection{Decomposing $S(0,N;n)$}  \label{sec:decomp_S2}

Here we will extract main term and error term behavior from
\begin{align*}
    S(0,N;n)&:=\sum_{\substack{k=1\\k\text{ even}}}^N\sum_{\substack{0\leq h<k\\(h,k)=1}}e^{-\frac{2\pi inh}{k}}\int_{-\vartheta_{h,k}'}^{\vartheta_{h,k}''}H_{\alpha,\beta}\left(e^{\frac{2\pi ih}{k}-\frac{2\pi z}{k}}\right)e^{\frac{2\pi nz}{k}}d\theta.
\end{align*}
Towards this goal, we apply the  transformation law (1) in Lemma \ref{lem:trans_law} for $H_{\alpha,\beta}$ and even $k$. In the process, we replace the resulting term 
\[H_{\alpha, \beta} \left( \exp \left( \dfrac{2 \pi i h'}{k} - \dfrac{2 \pi}{z k} \right) \right)=\sum_{j=0}^\infty a(\alpha,\beta;j) \exp \left( \dfrac{2 \pi i h'j}{k} - \dfrac{2 \pi j}{z k} \right)\] 
with its series expansion and re-express $z=kw$, yielding
\begin{align*}
    S(0,N;n)
    =&\sum_{j=0}^\infty\sum_{\substack{k=1\\k\text{ even}}}^N\sum_{\substack{0\leq h<k\\(h,k)=1}} A_{h,k}(\alpha,\beta,j;n)a(\alpha,\beta;j)k^{-\frac{\alpha+\beta}{2}}\\
    &\cdot \int_{-\vartheta_{h,k}'}^{\vartheta_{h,k}''}w^{-\frac{\alpha + \beta}{2}}\exp\left[w\left((\alpha+2\beta)\frac{\pi}{12}+2\pi n\right)+\frac{1}{w}\left(-(\alpha+2\beta)\frac{\pi}{12k^2}-\frac{2\pi j}{k^2}\right)\right]d\theta.
\end{align*}
The coefficient of $1/w$ inside the above exponential is positive if and only if $j<-(\alpha+2\beta)/24$.
Accordingly, we let
\begin{equation}
    S(0,N;n)
    =Q(0,N;n)+R(0,N;n)\label{eq:S(0)_split}
\end{equation}
where $Q(0,N;n)$ is the sum over those $j<-(\alpha+2\beta)/24$, and $R(0,N;n)$ consists of the remaining terms in $S(0,N;n)$. 
We will show later that $Q(0,N;n)$ yields mostly main term behavior, whereas $R(0,N;n)$ yields error term behavior. 
On this note, we observe that if $\alpha+2\beta\geq 0$, then $Q(0,N;n)$ is in fact an empty sum, which corresponds to the weak growth near $e^{2 \pi i h/k}$ for $k$ even discussed at the end of Section \ref{sec:trans_law}. 
Thus in our analysis of $Q(0,N;n)$ we may assume $\alpha+2\beta<0$.
Additionally, we will assume $n>-(\alpha+2\beta)/24$, which guarantees that the coefficient on $w$ inside the exponential is positive.

\subsection{Decomposing $Q(0,N;n)$}\label{sec:decomposing_Q(1)}
We now break up $Q(0,N;n)$ into three parts: one which will be extended to a Bessel function in the following subsection, as well as two error terms. 
Towards this goal, we divide the intervals of integration into three parts according to
$$-\vartheta_{h,k}' = - \frac{1}{k(k+k_1)} \leq - \frac{1}{k(N+k)} <  \frac{1}{k(N+k)}  \leq \frac{1}{k(k+k_2)} = \vartheta_{h,k}''. $$
Omitting the integrands, the split becomes
\begin{align}
Q(0,N;n)
    &=\sum_{j< -\frac{\alpha+2\beta}{24}}\sum_{\substack{k=2\\k\text{ even}}}^N\sum_{\substack{0\leq h<k\\(h,k)=1}} A_{h,k}(\alpha,\beta,j;n)a(\alpha,\beta;j)\nonumber\left[\int_{-\frac{1}{k(N+k)}}^{\frac{1}{k(N+k)}}\,
+
\int_{-\frac{1}{k(k_1+k)}}^{-\frac{1}{k(N+k)}}\,
+
\int_{\frac{1}{k(N+k)}}^{\frac{1}{k(k_2+k)}}\,\right]\nonumber\\
&=:Q_{0}(0,N;n)+Q_{1}(0,N;n)+Q_{2}(0,N;n)\label{eq:Q(0)_split_into_3}.
\end{align}

\subsection{Extending $Q_0(0,N;n)$\label{sec:Extend_Q_0}}
Here we will extend the path of integration of $Q_0(0,N;n)$ to obtain a modified Bessel function.
Implementing the variable transformation $w=N^{-2}-i\theta,$ we have
\begin{align*}
    Q_0(0,N;n)= & \sum_{j< -\frac{\alpha+2\beta}{24}}\sum_{\substack{k=2\\k\text{ even}}}^N\sum_{\substack{0\leq h<k\\(h,k)=1}} \frac{A_{h,k}(\alpha,\beta,j;n)a(\alpha,\beta,j)}{ik^{\frac{\alpha + \beta}{2}}} \cdot \int_{N^{-2}-\frac{i}{k(N+k)}}^{N^{-2}+\frac{i}{k(N+k)}}
    w^{-\frac{\alpha + \beta}{2}} \\
    & \cdot \exp\left[w\left((\alpha+2\beta)\frac{\pi}{12}+2\pi n\right)+\frac{1}{w}\left(-(\alpha+2\beta)\frac{\pi}{12k^2}-\frac{2\pi j}{k^2}\right)\right]dw.
\end{align*}
We will now extend the path of integration to the  rectangle $R$ with vertices
\[\pm N^{-2}\pm \frac{i}{k(N+k)}.\]
Omitting the integrands, we write
\begin{align}
Q_0(0,N;n)
    = 2 \pi
    \sum_{j< -\frac{\alpha+2\beta}{24}}
    &\sum_{\substack{k=2\\k\text{ even}}}^N
    \sum_{\substack{0\leq h<k\\(h,k)=1}}
    \frac{A_{h,k}(\alpha,\beta,j;n)a(\alpha,\beta,j)}{k^{\frac{\alpha + \beta}{2}}}\label{eq:Q0(0)_decomp}\\
    &\cdot\left\{\frac{1}{2\pi i}\int_R-\frac{1}{2\pi i}\left[\int_{+,+}^{-,+}+\int_{-,+}^{-,-}+\int_{-,-}^{+,-}\right]\right\}\nonumber.
\end{align}
We name $L^*(0,j,k;n):=1/2\pi i\int_R$, and we name the three remaining integrals $J_1(0,j,k,N;n),$ $J_2(0,j,k,N;n),$ and $J_3(0,j,k,N;n),$ respectively.
We will give an exact description of $L^*(0,j,k;n)$ as a Bessel function, and  bound each $J_i(0,j,k,N;n)$.

\subsection{Expressing $L^*(0,j,k;n)$ as a Bessel function}\label{sec:even_bessel}
In this subsection, we  express $L^*(0,j,k;n)$  { in terms of} a modified Bessel function of the first kind.
 
In the remaining subsections, we will bound each of the error terms that we have encountered, demonstrating that the only significant contribution from $S(0,N;n)$ comes from $L^*(0,j,k;n)$.
To $L^*(0,j,k;n)$ we apply the variable transformation
\[u=w\left((\alpha+2\beta)\frac{\pi}{12}+2\pi n\right),\]
which gives
\begin{align*}
L^*(0,j,k;n)=&\left[(\alpha+2\beta)\frac{\pi}{12}+2\pi n\right]^{\frac{\alpha+\beta}{2}-1}
\frac{1}{2\pi i}\int_{R'}u^{-\left(\frac{\alpha+\beta}{2}\right)} \\
 &\qquad \cdot \exp\left[u+\frac{1}{u}\left\{\left[-(\alpha+2\beta)\frac{\pi}{12k^2}-\frac{2\pi j}{k^2}\right]\left[(\alpha+2\beta)\frac{\pi}{12}+2\pi n\right]\right\}\right]du.
\end{align*}
It follows { from the integral representation of the $I$-Bessel function \eqref{eq:bes_int}} that
\begin{equation}\label{eq:L0Bessel}
L^*(0,j,k;n)=
\left[-(\alpha+2\beta)\frac{\pi}{12k^2}-\frac{2\pi j}{k^2}\right]^{\frac{1}{2}-\frac{\alpha+\beta}{4}}
\left[(\alpha+2\beta)\frac{\pi}{12}+2\pi n\right]^{-\frac{1}{2}+\frac{\alpha+\beta}{4}}
I_{v}(s_0(j,k))
\end{equation}
where
\begin{align*}
v:=1-\frac{\alpha+\beta}{2},\qquad s_0(j,k):=2\sqrt{\left[-(\alpha+2\beta)\frac{\pi}{12k^2}-\frac{2\pi j}{k^2}\right]\left[(\alpha+2\beta)\frac{\pi}{12}+2\pi n\right]}.
\end{align*}

\subsection{Bounding $R(0,N;n)$\label{sec:bounding_R}}

In this section we will bound $R(0,N;n)$, which we recall is given by
\begin{align*}
    R(0,N;n)
    =&\sum_{j \geq - \frac{\alpha + 2\beta}{24}}^\infty\sum_{\substack{k=1\\k\text{ even}}}^N\sum_{\substack{0\leq h<k\\(h,k)=1}} A_{h,k}(\alpha,\beta,j;n)a(\alpha,\beta;j)\\
    &\cdot \int_{-\vartheta_{h,k}'}^{\vartheta_{h,k}''}(kw)^{-\frac{\alpha + \beta}{2}}\exp\left[w\left((\alpha+2\beta)\frac{\pi}{12}+2\pi n\right)+\frac{1}{w}\left(-(\alpha+2\beta)\frac{\pi}{12k^2}-\frac{2\pi j}{k^2}\right)\right]d\theta.
\end{align*}
We will show that  $ R(0,N;n)= O(N^{-\delta}) $  for some  $\delta >0.$
Setting
$$E := E(w,\alpha, \beta,n,k,j) := \exp\left[w\left((\alpha+2\beta)\frac{\pi}{12}+2\pi n\right)+\frac{1}{w}\left(-(\alpha+2\beta)\frac{\pi}{12k^2}-\frac{2\pi j}{k^2}\right)\right]$$
and splitting the integral into three parts as in (\ref{eq:Q(0)_split_into_3}), we obtain
\begin{align} \nonumber R(0,N;n) = & \sum_{j \geq - \frac{\alpha + 2\beta}{24}}^\infty\sum_{\substack{k=1\\k\text{ even}}}^N\sum_{\substack{0\leq h<k\\(h,k)=1}} A_{h,k}(\alpha,\beta,j;n)a(\alpha,\beta;j) \\ \nonumber
& \cdot\bigg[ 
\int_{-\frac{1}{k(N+k)}}^{\frac{1}{k(N+k)}}(kw)^{-\frac{\alpha + \beta}{2}} E d \theta
 + \sum_{\ell=k_1+k}^{N+k-1} \int_{-\frac{1}{k\ell}}^{-\frac{1}{k(\ell+1)}}(kw)^{-\frac{\alpha + \beta}{2}} E d \theta 
+ \sum_{\ell=k_2+k}^{N+k-1}
\int_{\frac{1}{k\ell}}^{\frac{1}{k(\ell+1)}}(kw)^{-\frac{\alpha + \beta}{2}} E d \theta \bigg] \\
\label{eq:line_1}
= &  \sum_{j \geq - \frac{\alpha + 2\beta}{24}}^\infty a(\alpha,\beta;j) \sum_{\substack{k=1\\k\text{ even}}}^N\int_{-\frac{1}{k(N+k)}}^{\frac{1}{k(N+k)}}(kw)^{-\frac{\alpha + \beta}{2}} E d \theta \sum_{\substack{0\leq h<k\\(h,k)=1}} A_{h,k}(\alpha,\beta,j;n) 
\\ \label{eq:line_2}
& +\sum_{j \geq - \frac{\alpha + 2\beta}{24}}^\infty a(\alpha,\beta;j)\sum_{\substack{k=1\\k\text{ even}}}^N\sum_{\substack{0\leq h<k\\(h,k)=1}} A_{h,k}(\alpha,\beta,j;n) \sum_{\ell=k_1+k}^{N+k-1} \int_{-\frac{1}{k\ell}}^{-\frac{1}{k(\ell+1)}}(kw)^{-\frac{\alpha + \beta}{2}} E d \theta  
\\ \label{eq:line_3}
& +\sum_{j \geq - \frac{\alpha + 2\beta}{24}}^\infty a(\alpha,\beta;j)  \sum_{\substack{k=1\\k\text{ even}}}^N\sum_{\substack{0\leq h<k\\(h,k)=1}} A_{h,k}(\alpha,\beta,j;n)\sum_{\ell=k_2+k}^{N+k-1}
\int_{\frac{1}{k\ell}}^{\frac{1}{k(\ell+1)}}(kw)^{-\frac{\alpha + \beta}{2}} E d \theta.
\end{align}
Our goal is to bound the expressions (\ref{eq:line_1}), (\ref{eq:line_2}), and (\ref{eq:line_3}). We  begin by bounding $|E|$. 
As a reminder, we have $w= N^{-2} - i \theta$ and $|\theta| \leq1/Nk$, so that  $\mathrm{Re}(1/w) \geq k^2/2$
and
$ \mathrm{Re}(w) = N^{-2}.$ 
Therefore for $j \geq -(\alpha + 2\beta)/24$,
$$ |E| \leq \exp\left[N^{-2}\left((\alpha+2\beta)\frac{\pi}{12}+2\pi n\right)-(\alpha+2\beta)\frac{\pi}{24}-\pi j\right].$$
Thus for $N>0$, we have
\beq \label{eq:exp_bound}  \exp\left[N^{-2}\left((\alpha+2\beta)\frac{\pi}{12}+2\pi n\right)-(\alpha+2\beta)\frac{\pi}{24}\right] = O(1
).\eeq
Noting that the series 
$\sum_{j=0}^\infty a(\alpha, \beta;j)e^{-\pi j}$
defining $H_{\alpha, \beta}(i/2)$ converges absolutely, we have that for fixed $n$, (\ref{eq:line_1}) admits the estimate 
$$ O\left( \sum_{\substack{k=1\\k\text{ even}}}^N\int_{-\frac{1}{k(N+k)}}^{\frac{1}{k(N+k)}}(kw)^{-\frac{\alpha + \beta}{2}}  d \theta \sum_{\substack{0\leq h<k\\(h,k)=1}} A_{h,k}(\alpha,\beta,j;n)\right) = O\left(\sum_{\substack{k=1\\k\text{ even}}}^N  \frac{N^{\frac{\alpha + \beta}{2} -1}}{k} \sum_{\substack{0\leq h<k\\(h,k)=1}} A_{h,k}(\alpha,\beta,j;n)\right).$$
If $\alpha + \beta < 0$, we can use the trivial Kloosterman sum bound $O(k)$ to show that 
\beq O\left(\sum_{\substack{k=1\\k\text{ even}}}^N  \frac{N^{\frac{\alpha + \beta}{2} -1}}{k} \sum_{\substack{0\leq h<k\\(h,k)=1}} A_{h,k}(\alpha,\beta,j;n)\right) = O\left(N^{\frac{\alpha + \beta}{2}}\right). \label{eq:neg_weight}
\eeq
If $\alpha +\beta = 0$, we can make use of Lemma \ref{lem:Kloosterman} to show that for fixed $n$,
\beq O\left(\sum_{\substack{k=1\\k\text{ even}}}^N  \frac{N^{\frac{\alpha + \beta}{2} -1}}{k} \sum_{\substack{0\leq h<k\\(h,k)=1}} A_{h,k}(\alpha,\beta,j;n)\right) =O\left(\sum_{\substack{k=1\\k\text{ even}}}^N\frac{k^{-1/3 + \ep}}{N}\right) = O(N^{-1/3+ \ep}). \label{eq:zero_weight} \eeq

The summands (\ref{eq:line_2}) and (\ref{eq:line_3})  are handled in a similar way. The only important difference is that we must first switch the order of summation as in \cite[p.~507]{Rademacher} to obtain
\beq \label{eq:R_final} \sum_{j \geq - \frac{\alpha + 2\beta}{24}}^\infty a(\alpha,\beta;j)\sum_{\substack{k=1\\k\text{ even}}}^N  \sum_{\ell=N+1}^{N+k-1} \int_{-\frac{1}{k\ell}}^{-\frac{1}{k(\ell+1)}}(kw)^{-\frac{\alpha + \beta}{2}} E \sum_{\substack{0\leq h<k\\(h,k)=1\\N-k<k_2\leq l-k}} A_{h,k}(\alpha,\beta,j;n) d \theta. 
\eeq
The desired bound is now obtained via the same methods used to bound (\ref{eq:line_1}). As described in Section \ref{sec:Kloosterman},  Lemma \ref{lem:Kloosterman} is equipped to handle the incomplete Kloosterman sum in (\ref{eq:R_final}). By (\ref{eq:neg_weight}), (\ref{eq:zero_weight}), and corresponding statements for (\ref{eq:line_2}) and (\ref{eq:line_3}), we obtain in all cases
\beq \label{eq:R0_bound} R(0,N;n) = O(N^{-\delta}) \eeq
for some $\delta >0$.

\subsection{Bounding $Q_1(0,N;n)$ and $Q_2(0,N;n)$}\label{sec:boundingQ_2(0)}

In this subsection, we will bound those segments of the Farey arcs that do not contribute to the Bessel integral.
We will explicitly bound just $Q_2(0,N;n)$ because similar arithmetic yields the same bound for $Q_1(0,N;n)$.
First, as in (\ref{eq:line_2}) and (\ref{eq:line_3}) we split our path of integration into many smaller intervals, obtaining
\begin{align*}
Q_2(0,N;n)
    =&\sum_{j< -\frac{\alpha+2\beta}{24}}^\infty
    \sum_{\substack{k=1\\k\text{ even}}}^N
    \sum_{\substack{0\leq h<k\\(h,k)=1}} A_{h,k}(\alpha,\beta,j;n)
    a(\alpha,\beta;j)k^{-\frac{\alpha + \beta}{2}}
    \sum_{l=k_2+k}^{N+k-1} \int_{\frac{1}{k(l+1)}}^{\frac{1}{kl}}
    w^{-\frac{\alpha + \beta}{2}}
    E d\theta.
\end{align*}
Now, as in \ref{eq:R_final}, we switch the order of summation, which yields
\begin{align*}
Q_2(0,N;n)
    =&\sum_{j< -\frac{\alpha+2\beta}{24}}^\infty a(\alpha,\beta;j)
    \sum_{\substack{k=1\\k\text{ even}}}^N
    k^{-\frac{\alpha + \beta}{2}}
    \sum_{l=N+1}^{N+k-1} 
    \int_{\frac{1}{k(l+1)}}^{\frac{1}{kl}}
    w^{-\frac{\alpha + \beta}{2}}
    Ed\theta \sum_{\substack{0\leq h<k\\(h,k)=1\\N-k<k_2\leq l-k}}
    A_{h,k}(\alpha,\beta,j;n).
\end{align*}
We have, similarly to (\ref{eq:exp_bound}), that $E=O(1)$. 
Since 
$ |w| \leq ({2}/{(N^2k^2)})^{1/2}$,
we have
\begin{align*}
    \sum_{l=N+1}^{N+k-1}\int_{\frac{1}{k(l+1)}}^{\frac{1}{kl}}
    &w^{-\frac{\alpha + \beta}{2}}
    Ed\theta = O\left(\left(Nk\right)^{\frac{\alpha+\beta}{2}-1}\right).\\
\end{align*}
Here and throughout, the constant implied by the big Oh notation depends at most on $\alpha,\beta,n$.
The following steps are analogous to the final steps of Section \ref{sec:bounding_R}, where the cases $\alpha + \beta <0$ and $\alpha + \beta =0$ are considered separately. 
We conclude that
\begin{equation}
    Q_1(0,N;n), \ Q_2(0,N;n)=O(N^{-\delta})\label{eq:Q_2(0)_bound}
\end{equation}
for some $\delta>0$.

\subsection{Bounding $J_1(0,j,k,N;n),J_2(0,j,k,N;n),J_3(0,j,k,N;n)$\label{sec:boundingJs}}

As a reminder, we have 
\[J_1(0,j,k,N;n):=\int_{N^{-2}+\frac{i}{k(N+k)}}^{-N^{-2}+\frac{i}{k(N+k)}}w^{-\frac{\alpha + \beta}{2}}
    Edw, \qquad J_3(0,j,k,N;n):=\int_{-N^{-2}-\frac{i}{k(N+k)}}^{+N^{-2}-\frac{i}{k(N+k)}}w^{-\frac{\alpha + \beta}{2}}
   Edw,\]
and 
\[J_2(0,j,k,N;n):=\int_{-N^{-2}+\frac{i}{k(N+k)}}^{-N^{-2}-\frac{i}{k(N+k)}}w^{-\frac{\alpha + \beta}{2}}
    Edw.\]
    
We first consider $J_1$ and $J_3$. Here 
one easily finds $\mathrm{Re}(w) \leq N^{-2}$ and $ \mathrm{Re}(1/w) \leq 4k^2$, and from this obtains $E= O(1)$.
Since the interval length is $2N^{-2}$ and $|w| \leq {\sqrt{2}}/{kN}$,
 we have that $J_1$ and $J_3$ are $O\left( (kN)^{(\alpha + \beta)/2-1})\right).$

For $J_2$, we need only note that $\mathrm{Re}(w)$ and $ \mathrm{Re}(1/w) $ are negative to show that $E= O(1)$. Since the length of the interval is $2/N(k+N)$, we see that  $J_2$ is $O\left( (kN)^{(\alpha + \beta)/2-1})\right).$
We proceed as in Sections \ref{sec:bounding_R} and \ref{sec:boundingQ_2(0)} to show that 
\beq \label{eq:J_bound} \sum_{\substack{k=2 \\ k \ \text{even}}}^NA_k(\alpha,\beta,j;
n)\left\{J_1+J_2+J_3\right\} = O(N^{-\delta}) \eeq
for some $\delta>0$, as desired.

Recalling our decomposition of $Q_0(0,N;n)$ in  (\ref{eq:Q0(0)_decomp}) and the bound given by (\ref{eq:J_bound}) above, we can now estimate
\begin{align}
    Q_0(0,N;n)=
    2\pi\sum_{j< -\frac{\alpha+2\beta}{24}}
    \sum_{\substack{k=2\\k\text{ even}}}^N
    \frac{A_{k}(\alpha,\beta,j;n)a(\alpha,\beta,j)}{k^{\frac{\alpha + \beta}{2}}}L^*(0,j,k;n)+O(N^{-\delta}).\label{eq:Q0(0)_estimated}
\end{align}

\section{The odd case} \label{sec:odd_case}

The odd case follows a manner very similar to the even case, and the discrepancies are all consequences of the differences between the transformation formulas (1) and (2)  in Lemma \ref{lem:trans_law} for even and odd $k$.

\subsection{Decomposing $S(1,N;n)$}

Here we will extract main term and error term behavior from
\begin{align*}
    S(1,N;n)&=\sum_{\substack{0\leq h<k\leq N\\(h,k)=1\\k\text{ odd}}}e^{-\frac{2\pi inh}{k}}\int_{-\vartheta_{h,k}'}^{\vartheta_{h,k}''}H_{\alpha,\beta}\left(e^{\frac{2\pi ih}{k}-\frac{2\pi z}{k}}\right)e^{\frac{2\pi nz}{k}}d\theta.
\end{align*}
We apply the  transformation law (2) in Lemma \ref{lem:trans_law} for $H_{\alpha,\beta}$ and odd $k$, obtaining 
\begin{align*}
    S(1,N;n)=&
    \sum_{j=0}^\infty 
    \sum_{\substack{k=1\\k\text{ odd}}}^N
    \sum_{\substack{0\leq h<k\\(h,k)=1}} \frac{B_{h,k}(\alpha,\beta,j;n)a(\beta,\alpha;j)}{2^{\beta/2}k^{\frac{\alpha+\beta}{2}}}\\
    &\cdot \int_{-\vartheta_{h,k}'}^{\vartheta_{h,k}''}w^{-\frac{\alpha + \beta}{2}}
    \exp\left[w\left((\alpha+2\beta)\frac{\pi}{12}+2\pi n\right)
    +\frac{1}{w}\left(-(2\alpha+\beta)\frac{\pi}{24k^2}-\frac{\pi j}{k^2}\right)\right]d\theta
\end{align*}
as in Section \ref{sec:decomp_S2}.

In this case the coefficient on $1/w$ inside the above exponential is positive if and only if $j<-(2\alpha+\beta)/24.$
Thus we stratify $ S(1,N;n)$ via
\begin{equation}
    S(1,N;n)
    =Q(1,N;n)+R(1,N;n),\label{eq:S(1)_split}
\end{equation}
where $Q(1,N;n)$ is the sum over those $j<-(2\alpha+\beta)/24$.
As in the even case,  $Q(0,N;n)$ will yield mostly main term behavior, whereas $R(0,N;n)$ will yield error term behavior. 
We continue to assume that $n>-(\alpha+2\beta)/24$.

\subsection{Decomposing $Q(1,N;n)$}

We will now break up $Q(1,N;n)$ into three parts, which play the same roles as in Section \ref{sec:decomposing_Q(1)}. 
Omitting the integrands, we have
\begin{align}
    Q(1,N;n)&=
    \sum_{j<-\frac{2\alpha+\beta}{24}}
    \sum_{\substack{k=1\\k\text{ odd}}}^N
    \sum_{\substack{0\leq h<k\\(h,k)=1}} \frac{B_{h,k}(\alpha,\beta,j;n)a(\beta,\alpha;j)}{2^{\beta/2}k^{\frac{\alpha+\beta}{2}}}\nonumber\cdot\left[\int_{-\frac{1}{k(N+k)}}^{\frac{1}{k(N+k)}}\,
+
\int_{-\frac{1}{k(k_1+k)}}^{-\frac{1}{k(N+k)}}\,
+
\int_{\frac{1}{k(N+k)}}^{\frac{1}{k(k_2+k)}}\,\right]\\
&=:Q_{0}(1,N;n)+Q_{1}(1,N;n)+Q_{2}(1,N;n).\label{eq:Q(1)_split_into_3}
\end{align}

\subsection{Extending $Q_0(1,N;n)$}

In this subsection, we will extend the path of integration of $Q_0(1,N;n)$ to obtain a modified Bessel function.
Via the variable transformation $w=N^{-2}-i\theta$, we can write
\begin{align*}
    Q_0(1,N;n)=&
    \sum_{j<-\frac{2\alpha+\beta}{24}}
    \sum_{\substack{k=1\\k\text{ odd}}}^N
    \sum_{\substack{0\leq h<k\\(h,k)=1}} \frac{B_{h,k}(\alpha,\beta,j;n)a(\beta,\alpha;j)}{i2^{\beta/2}k^{\frac{\alpha+\beta}{2}}}
    \\
    &\cdot \int_{N^{-2}-\frac{i}{k(N+k)}}^{N^{-2}+\frac{i}{k(N+k)}}
    w^{-\frac{\alpha + \beta}{2}}
    \exp\left[w\left((\alpha+2\beta)\frac{\pi}{12}+2\pi n\right)
    +\frac{1}{w}\left(-(2\alpha+\beta)\frac{\pi}{24k^2}-\frac{\pi j}{k^2}\right)\right]dw.
\end{align*}
Considering the rectangle $R$ described in Section \ref{sec:Extend_Q_0},
we write
\begin{align}
    Q_0(1,N;n)=&
    2 \pi\sum_{j<-\frac{2\alpha+\beta}{24}}
    \sum_{\substack{k=1\\k\text{ odd}}}^N
    \sum_{\substack{0\leq h<k\\(h,k)=1}} \frac{B_{h,k}(\alpha,\beta,j;n)a(\beta,\alpha;j)}{2^{\beta/2}k^{\frac{\alpha+\beta}{2}}}\label{eq:Q0(1)_decomp}\cdot \left\{\frac{1}{2\pi i}\int_R-\frac{1}{i}\left[\int_{+,+}^{-,+}+\int_{-,+}^{-,-}+\int_{-,-}^{+,-}\right]\right\}.
\end{align}
Omitting the integrands, we write $L^*(1,j,k;n):=1/2\pi i\int_R$, and we name the three remaining integrals $J_1(1,j,k,N;n),$ $J_2(1,j,k,N;n),$ and $J_3(1,j,k,N;n),$ respectively.

\subsection{Expressing $L^*(1,j,k;n)$ as a Bessel function}

Now, we will express $L^*(1,j,k;n)$ exactly as a modified Bessel function of the first kind.
To $L^*(1,j,k;n)$ we apply the same variable transformation as in Section \ref{sec:even_bessel}, which gives
\begin{equation}\label{eq:L1Bessel}
    L^*(1,j,k;n)
    =
    \left[-(2\alpha+\beta)\frac{\pi}{24k^2}-\frac{\pi j}{k^2}\right]^{\frac{1}{2}-\frac{\alpha+\beta}{4}}\left[(\alpha+2\beta)\frac{\pi}{12}+2\pi n\right]^{-\frac{1}{2}+\frac{\alpha+\beta}{4}}I_v(s_1(j,k)),
\end{equation}
where 
\[s_1(j,k):=2\sqrt{\left[-(2\alpha+\beta)\frac{\pi}{24k^2}-\frac{\pi j}{k^2}\right]\left[(\alpha+2\beta)\frac{\pi}{12}+2\pi n\right]}\]
and $v$ is the same as in Section \ref{sec:even_bessel}.

\subsection{Bounding the odd case error terms}

 $R(1,N;n)$, $Q_\mu(1,N;n)$, and $J_\nu(1,j,k,N;n)$ for $\mu=1,2$ and $\nu=1,2,3$ are bounded in a manner which is nearly identical to the corresponding arguments in the even case. The main differences between these expressions and their even case counterparts are the coefficient $a(\beta,\alpha;j)$ from the expansion of $H_{\beta,\alpha}$ and the coefficient 
 $$-(2\alpha +\beta)\frac{\pi}{24k^2}-\frac{\pi j}{k^2}$$
 of $1/w$ in the exponential. Note that in Sections \ref{sec:bounding_R}, \ref{sec:boundingQ_2(0)}, and \ref{sec:boundingJs}, the terms $-(\alpha + 2\beta)/12$ in the coefficient of $1/w$, which correspond to $-(2\alpha +\beta)/24$ in the odd case, may as well have been arbitrary negative real numbers that are fixed in this discussion. The same holds for $-2 \pi$ as the coefficient of $j$, except when we used the absolute convergence of the series defining $H_{\alpha,\beta}(i/2)$ in Section \ref{sec:bounding_R}. 
 In the odd case, we need to use the absolute convergence of the series defining $H_{\beta,\alpha}(i/4)$ to obtain the same result. Running the same arguments with these minor adjustments, we obtain
 \begin{align}
    R(1,N;n) & = O(N^{-\delta}),  \label{eq:R(1)_bound} \\
    Q_\mu(1,N;n) & = O(N^{-\delta}),  \label{eq:Q_2(1)_bound}  \end{align}
    and
    \begin{align}
    \sum_{\substack{0\leq h<k\\(h,k)=1}} B_{h,k}(\alpha,\beta,j;n) \left( \sum_\nu J_\nu(1,j,k,N;n)\right) & = O(N^{-\delta})  \label{eq:J(1)_bound}
 \end{align}
 for some $\delta>0$, where $ \mu=1,2$ and $\nu=1,2,3$.

Recalling our decomposition of $Q_0(1,N;n)$ in (\ref{eq:Q0(1)_decomp}) and the bounds given by (\ref{eq:J(1)_bound}) above, we conclude that
\begin{align}
    Q_0(1,N;n)=
    2\pi\sum_{j< -\frac{2\alpha+\beta}{24}}
    \sum_{\substack{k=1\\k\text{ odd}}}^N
    \frac{B_{k}(\alpha,\beta,j;n)a(\beta,\alpha,j)}{2^{\beta/2}k^{\frac{\alpha + \beta}{2}}}L^*(1,j,k;n)+O(N^{-\delta}).\label{eq:Q0(1)_estimated}
\end{align}

\section{Proof of Theorem \ref{thm:main_result} and Corollaries \ref{cor:asymptotics} and \ref{cor:equidistribution}} \label{sec:proofs}

We now apply our work from Sections \ref{sec:Circle_method_outline}, \ref{sec:Kloosterman}, and \ref{sec:even_case}  to obtain exact formulas for $a(\alpha,\beta;n)$, and afterwards we will extract asymptotic behavior from those formulas.
Finally, we will express these exact and asymptotic formulas in terms of our topological invariants to prove Theorem \ref{thm:main_result} and Corollary \ref{cor:asymptotics}.

\begin{theorem}\label{thm:exactFormulasFor_H}
Let $\alpha+\beta\leq 0$, and assume $n>-(\alpha+2\beta)/24$.
Then
\begin{align} \label{eq:exactFormulasFor_H}
    a(\alpha,\beta;n)=\ \ \ &2\pi
    \sum_{j<-\frac{\alpha+2\beta}{24}}
    \sum_{\substack{k=2\\k\mathrm{\ even}}}^\infty
    \frac{A_{k}(\alpha,\beta,j;n)a(\alpha,\beta,j)}{k^{\frac{\alpha + \beta}{2}}}L^*(0,j,k;n)\\
    \qquad+\,&2\pi\sum_{j<-\frac{2\alpha+\beta}{24}}
    \sum_{\substack{k=1\\k \ \mathrm{ odd}}}^\infty \frac{B_{k}(\alpha,\beta,j;n)a(\beta,\alpha;j)}{2^{\beta/2}k^{\frac{\alpha+\beta}{2}}}L^*(1,j,k;n).
\end{align}
\end{theorem}
\begin{proof}
Making use of the decompositions (\ref{eq:S(0)_split}) and (\ref{eq:Q(0)_split_into_3}), as well as the estimates (\ref{eq:R0_bound}), (\ref{eq:Q_2(0)_bound}),  and (\ref{eq:Q0(0)_estimated}), we have 
\begin{align*}
    S(0,N;n)=2\pi\sum_{j< -\frac{\alpha+2\beta}{24}}
    \sum_{\substack{k=2\\k\text{ even}}}^N
    \frac{A_{k}(\alpha,\beta,j;n)a(\alpha,\beta,j)}{k^{\frac{\alpha + \beta}{2}}}L^*(0,j,k;n)+O(N^{-\delta}).
\end{align*}
Similarly, by the decompositions (\ref{eq:S(1)_split}) and (\ref{eq:Q(1)_split_into_3}), as well as the estimates (\ref{eq:R(1)_bound}), (\ref{eq:Q_2(1)_bound}), and (\ref{eq:Q0(1)_estimated}), we have 
\begin{align*}
    S(1,N;n)=
    2\pi\sum_{j< -\frac{2\alpha+\beta}{24}}
    \sum_{\substack{k=1\\k\text{ odd}}}^N
    \frac{B_{k}(\alpha,\beta,j;n)a(\beta,\alpha,j)}{2^{\beta/2}k^{\frac{\alpha + \beta}{2}}}L^*(1,j,k;n)+O(N^{-\delta})\label{eq:Q0(1)_estimated}.
\end{align*}
But by the parity split in (\ref{eq:a_parity_strat_S's}), these estimates imply that
\begin{align*}
    a(\alpha,\beta;n)=\ \ \ &2\pi
    \sum_{j<-\frac{\alpha+2\beta}{24}}
    \sum_{\substack{k=2\\k\mathrm{\ even}}}^N
    \frac{A_{k}(\alpha,\beta,j;n)a(\alpha,\beta,j)}{k^{\frac{\alpha + \beta}{2}}}L^*(0,j,k;n)\\
    \qquad+\,&2\pi\sum_{j<-\frac{2\alpha+\beta}{24}}
    \sum_{\substack{k=1\\k \ \mathrm{ odd}}}^N \frac{B_{k}(\alpha,\beta,j;n)a(\beta,\alpha;j)}{2^{\beta/2}k^{\frac{\alpha+\beta}{2}}}L^*(1,j,k;n)\\
    \qquad+\,&O(N^{-\delta})
\end{align*}
for some $\delta>0$.
Keeping $n$ fixed and letting $N\to\infty$, we obtain (\ref{eq:exactFormulasFor_H}).
\end{proof}

{ \noindent \tbf{Remark:}
In \cite{BringOno}, Bringmann and Ono obtained exact formulas for the coefficients of harmonic Maass forms of non-positive weight using the theory of Maass-Poincar\'e series.
Our results can therefore also be obtained using their work. \medskip
}

In the following corollary, we make use of Theorem \ref{thm:exactFormulasFor_H} to obtain approximations for $a(\alpha,\beta;n)$ for large $n$. We will see that the main contribution for large $n$ comes from one of the first two summands in $k$. 

\begin{corollary}\label{cor:asym_for_h} 
Suppose that $\alpha + \beta \leq 0$. 
\begin{enumerate}
\item If $\alpha = 0$ and $\beta \neq 0$, then as $n \to \infty$ we have
 $$
 a(0, \beta;2n) \sim 2^{\frac{3\beta-5}{4}}3^{\frac{\beta-1}{4}}(-\beta)^{\frac{1-\beta}{4}}n^{\frac{\beta-3}{4}}\exp\left(\pi\sqrt{\frac{-2\beta}{3}n}\right).$$
\item If
$\alpha<0$,  then as $n \to \infty$ we have
$$ a(\alpha,\beta;n) \sim 2^{\frac{2\alpha+\beta-3}{2}}3^{\frac{\alpha+\beta-1}{4}}(-(2\alpha+\beta))^{\frac{-\alpha-\beta+1}{4}}n^{\frac{\alpha+\beta-3}{4}}\exp\left(\pi\sqrt{\frac{-(2\alpha+\beta)}{3}n}\right).$$
\item If 
$\alpha>0$, then as $n \to \infty$ we have
$$ a(\alpha, \beta;n) \sim (-1)^{n} 2^{\frac{3\alpha+3\beta-7}{4}}3^{\frac{\alpha+\beta-1}{4}}(-(\alpha+2\beta))^{\frac{1-\alpha-\beta}{4}}n^{\frac{\alpha+\beta-3}{4}}\exp\left(\pi\sqrt{\frac{-(\alpha+2\beta)}{6}n}\right) .$$
\end{enumerate}
\end{corollary}
\begin{proof}
Consider the case $\alpha < 0.$ Let
\begin{equation*}
    a(\alpha, \beta; n) =: 2 \pi \bigg( M^-(\alpha, \beta; n) + E^-(\alpha, \beta; n) \bigg),
\end{equation*}
where we let
\begin{equation*}
    M^-(\alpha, \beta; n) := \dfrac{B_1(\alpha, \beta, 0; n) a(\beta, \alpha; 0)}{2^{\frac{\beta}{2}}} L^*(1, 0, 1; n)
\end{equation*}
and
\begin{align}\label{eq:weak_bess}
    E^-(\alpha, \beta;  n) := 2 \pi \left(\sum_{j = 0}^{\mathlarger{\lfloor} \frac{\alpha+2\beta}{24} \mathlarger{\rfloor}}
    \sum_{\substack{k=2\\k \equiv 0 \Mod{2}}}^\infty f_A(\alpha, \beta, j, k; n) + \sum_{j = 0} ^{\mathlarger{\lfloor} \frac{2\alpha+\beta}{24} \mathlarger{\rfloor}}
    \sum_{\substack{k=3\\k \equiv 1 \Mod{2}}}^\infty f_B(\alpha, \beta, j, k; n) \right),
\end{align}
where
\begin{equation*}
    f_A(\alpha, \beta, j, k; n) := \dfrac{A_k(\alpha, \beta, j; n) a(\alpha, \beta, j)}{k^{\frac{\alpha + \beta}{2}}} L^*(0, j, k; n)
\end{equation*}
and
\begin{equation*}
    f_B(\alpha, \beta, j, k; n) := \dfrac{B_k(\alpha, \beta, j; n) a(\beta, \alpha, j)}{2^\beta k^{\frac{\alpha + \beta}{2}}} L^*(1, j, k; n).
\end{equation*}
To prove Corollary \ref{cor:asym_for_h} for $\alpha < 0,$ it suffices to show that $E^-(\alpha, \beta; n) = o(L^*(1,0,1;n))$ and then to carry out the necessary simplification of $2 \pi M^-(\alpha, \beta; n).$

{ Using the bound \eqref{eq:geom_bes_bound}}, the monotonicity of the $I$-Bessel functions, and the trivial bound on the Kloosterman sums, it follows from \eqref{eq:weak_bess} that
\begin{align*}
    E^-(\alpha, \beta; n) &\leq \bigg( \dfrac{-(\alpha+2\beta)}{\alpha+2\beta+24 n} \bigg)^{\frac{v}{2}} \sum_{j = 0}^{\mathlarger{\lfloor} -\frac{\alpha + 2 \beta}{24} \mathlarger{\rfloor}} |a(\alpha,\beta,j)| \bigg( K_A(\alpha, \beta, j; n) + \lceil s_0(j, 1) \rceil I_{v}(s_0(0,2) ) \bigg)
    \\
    & + \bigg( \dfrac{-(2\alpha+\beta)}{\alpha+2\beta + 48 \pi n} \bigg)^{\frac{v}{2}} \sum_{j = 0}^{\mathlarger{\lfloor} -\frac{2 \alpha + \beta}{24} \mathlarger{\rfloor}} \frac{ |a(\beta,\alpha;j)| }{2^{\beta/2}} \bigg( K_B(\alpha, \beta, j; n) + \lceil s_1(j,1) \rceil I_{v}( s_1(0,3)) \bigg),
\end{align*}
where 
\begin{equation*}
    K_A(\alpha, \beta, j; n) := \dfrac{4}{3} \sum_{\substack{k > \lceil s_0(j, 1) \rceil \\ k \equiv 0 \Mod{2}}}\frac{A_{k}(\alpha,\beta,j;n)}{k^{v + 1}} \bigg( \dfrac{s_0(0,1)}{2} \bigg)^{v},
\end{equation*}
and
\begin{equation*}
    K_B(\alpha, \beta, j; n) := \dfrac{4}{3} \sum_{\substack{k > \lceil s_1(j,1) \rceil \\ k \equiv 1 \Mod{2}}} \frac{B_{k}(\alpha,\beta,j;n)}{k^{v + 1}} \bigg( \dfrac{s_1(0,1)}{2} \bigg)^{v}.
\end{equation*}
If $\alpha + \beta <0$, we use the trivial bound on the Kloosterman sums to obtain
\begin{equation*}
    K_A(\alpha, \beta, j; n) \leq \dfrac{4}{3}  \bigg( \dfrac{s_0(0,1)}{2} \bigg)^{v} \dfrac{s_0(j,1)^{1 - v}}{1 - v}
\end{equation*}
and
\begin{equation*}
    K_B(\alpha, \beta, j; n) \leq \dfrac{4}{3} \bigg( \dfrac{s_1(0,1)}{2} \bigg)^{v} \dfrac{s_0(j,1)^{1 - v}}{1 - v},
\end{equation*}
and if $\alpha + \beta = 0,$ we use Lemma \ref{lem:Kloosterman} to obtain
\begin{equation*}
    K_A(\alpha, \beta, j; n) =  O(n^{\frac{1}{3}} s_0(j, 1)^{-\frac{1}{3} + \epsilon})
\end{equation*}
and
\begin{equation*}
    K_B(\alpha, \beta, j; n) =  O(n^{\frac{1}{3}} s_1(j, 1)^{-\frac{1}{3} + \epsilon}).
\end{equation*}
Therefore, since $E^-(\alpha, \beta; n)$ is bounded by a finite sum of rational functions and $I$-Bessel functions, it follows from
\begin{equation*}
    I_v(x)\sim \frac{e^x}{\sqrt{ 2\pi x}}
\end{equation*}
that
$E^-(\alpha, \beta; n) = o(L^*(1,0,1;n)),$ and that the statement of Corollary \ref{cor:asym_for_h} follows from a simplification of $2 \pi M(\alpha, \beta; n).$

If $\alpha > 0,$ the proof of Corollary \ref{cor:asym_for_h} follows \emph{mutatis mutandis} if we let 
\begin{equation*}
    a(\alpha, \beta; n) := 2 \pi \bigg( M^+(\alpha, \beta; n) + E^+(\alpha, \beta; n) \bigg),
\end{equation*}
where we let
\begin{equation*}
    M^+(\alpha, \beta; n) := \dfrac{A_2(\alpha, \beta, 0; n) a(\beta, \alpha; 0)}{2^{\frac{\alpha +  \beta}{2}}} L^*(0, 0, 2; n)
\end{equation*}
and
\begin{align*}
    E^+(\alpha, \beta;  n) := 2 \pi \left(\sum_{j = 0}^{\mathlarger{\lfloor} \frac{\alpha+2\beta}{24} \mathlarger{\rfloor}}
    \sum_{\substack{k=4\\k \equiv 0 \Mod{2}}}^\infty f_A(\alpha, \beta, j, k; n) + \sum_{j = 0} ^{\mathlarger{\lfloor} \frac{2\alpha+\beta}{24} \mathlarger{\rfloor}}
    \sum_{\substack{k=1\\k \equiv 1 \Mod{2}}}^\infty f_B(\alpha, \beta, j, k; n) \right).
\end{align*}
Finally, we observe that $a(0, \beta; 2 n)  =  a(\beta, 0, n),$ and so making the necessary substitutions for the case $\alpha < 0$ yields the  case $\alpha = 0.$
\end{proof}

We are now in a position to prove Theorem \ref{thm:main_result}, our exact formulas for $\chi(\mathrm{Hilb}^n(S))$ and $\sigma(\mathrm{Hilb}^n(S))$, and Corollary \ref{cor:asymptotics}, our asymptotic formulas for $\chi(\mathrm{Hilb}^n(S))$, $\sigma(\mathrm{Hilb}^n(S))$, $b^*(r,2;n)$ and $c^*(r_1,r_2;n)$.
At this point, the work is a simple application of Theorem \ref{thm:exactFormulasFor_H} and Corollary \ref{cor:asym_for_h}.

\begin{proof}[Proof of Theorem \ref{thm:main_result}]
The proof follows from Lemma \ref{lem:gottsche_prod} and Theorem \ref{thm:exactFormulasFor_H}. We note only that in the derivation of (\ref{eq:Euler_char_exact_formula}) the sum over odd $k$ is actually a sum over even $j$, since $a(0,-\chi(S);j) =0$ for odd $j$. Noting this and replacing $j$ with $2j$, we obtain  (\ref{eq:Euler_char_exact_formula}). 
\end{proof}


\begin{proof}[Proof of Corollary \ref{cor:asymptotics}]
This follows from Corollary
\ref{cor:asym_for_h} and 
(\ref{sp2}) 
and (\ref{sp3}) in Lemma \ref{lem:gottsche_prod}. 
\end{proof}

\begin{proof}[Proof of Corollary \ref{cor:equidistribution}]
The coefficients of our functions $Z_S(x, y;\tau)$, where $x,y = \pm1$, are of the form $a(\alpha,\beta;n)$ for some $\alpha$ and $\beta$ determined by (\ref{eq:minus_expression}) and Lemma \ref{lem:gottsche_prod}. For convenience we define 
$$a_S(x, y;n):= a(\alpha,\beta;n).$$ By Corollary
\ref{cor:asym_for_h}, each of these coefficient sequences is asymptotic to a function of the form $c_1n^{c_2}\exp(\sqrt{Gn})$, where $c_1$, $c_2$, and $G$ only depend $x$, $y$, and $S$. We will let $G_S(x ,y)$ denote this value $G$. Note that $G_S(x,y)$ essentially determines the growth of  the sequence $a_S(x,y;n)$ in the sense that if $G_S(x_1,y_1)<G_S(x_2,y_2)$, then $a_S(x_1, y_1;n) = o(a_S(x_2,y_2;n))$. 

The proof of (1) now follows easily from Corollary
\ref{cor:asym_for_h}, 
(\ref{sp2}), (\ref{sp3}), and (\ref{eq:minus_expression}).
As a reminder, (\ref{eq:B_sum}) tells us 
$$ B_S(r,2;\tau)= \frac{1}{2}\left( Z(1,1;\tau) + (-1)^r Z(1,-1;\tau)\right) .
$$
One makes use of Corollary \ref{cor:asym_for_h} to check that $G_S,(1,1)>G_S(1,-1)$  if and only if $\sigma(S) + \chi(S)>0$ and $G_S,(1,1)<G_S(1,-1)$ if and only if $\sigma(S) + \chi(S)<0$. Note that $G_S(1,1)=G_S(1,-1)$ whenever $\chi(S) = - \sigma(S)$, and in this case  $Z_S(1,1;\tau)=Z_S(1,-1;\tau)$, so that  $b^*_S(1;n)=0$ for all $n$. Also,   if $\chi(S) \leq 0$, $(\chi(S) ,\sigma(S)) \neq (0,0)$,  we have from Theorem 15.1 in \cite{HMF} that $a_S(1,1;n) = o(a_S(1,-1;n))$. Since this requires $\sigma(S) + \chi(S)\leq 0$, we can remove the hypothesis that $\chi(S)> 0$. 

To prove (2), we first recall that (\ref{eq:now}) gives
$$C_S(r_1, r_2;\tau) = \frac{1}{4} \sum_{\substack{ j_1 \Mod{2} \\ j_2\Mod{2}}}  (-1)^{j_2r_2} (-1)^{j_1r_1} Z_S((-1)^{j_2}, (-1)^{j_1};\tau).$$
We note that by (\ref{eq:minus_expression}), the hypotheses of Corollary \ref{cor:asym_for_h} are always satisfied for $Z(-1,-1;\tau)$, since 
$$4h^{1,0} - (\chi(S) -8h^{1,0})= -(2h^{0,0} + 2h^{2,0} + h^{1,0})\leq0$$
 for all $S$. Note also that for $\chi(S) \geq 0$, we have
 $$G_S(-1,-1) - G_S(1,1) = \frac{2\chi(S) + 12h^{1,0}}{3} - \frac{2\chi(S)}{3} = 4h^{1,0} \geq 0.$$
  If $h^{1,0}=0$, we have $Z_S(1,1;\tau)=Z_S(-1,-1;\tau)$. If $\chi(S) \leq 0$,  we have     $a_S(1,1;n) = o(a_S(-1,-1;n))$ by Theorem 15.1 in \cite{HMF}. Since $\chi(S)\leq 0$ requires $h^{1,0}>0$, we have  
  $$a_S(1,1;n) + (-1)^{r_1+r_2}a_S(-1,-1;n) \sim K a_S(-1,-1;n), $$ where 
  $$ K : = 
  \bcas \hspace{2.9mm}  2 & h^{1,0} = 0  \text{ and } r_1 + r_2 \equiv 0 \mod 2 \\
  \hspace{2.9mm} 0 & h^{1,0}=0 \text{ and } r_1 + r_2 \equiv 1 \mod 2 \\
  \hspace{3.1mm} 1 & h^{1,0} > 0 \text{ and } r_1 + r_2 \equiv 0 \mod 2 \\
  -1 & h^{1,0} > 0 \text{ and } r_1 + r_2 \equiv 1 \mod 2. \\
  \ecas
  $$
We now see that if 
$G_S(-1,-1)> G_S(1,-1) $, then  $c^*_S(r_1,r_2;n) \sim (K/4)a_S(-1,-1;n)$, which gives the desired result. 
To conclude the proof, one simply checks that this inequality holds for both $\sigma(S) \geq 0$ and $\sigma(S)<0$.   
\end{proof}


\section{Examples}\label{sec:numerics}

Here we illustrate Theorem \ref{thm:main_result} 
and Corollaries \ref{cor:asymptotics} and \ref{cor:equidistribution} 
with examples of numerical computations.

\begin{example} 
\emph{ To demonstrate Theorem \ref{thm:main_result}, we focus on $S$, a 2-dimensional torus blown up at one point, a Hirzebruch surface $\Sigma_0 \cong \mathbb{P}^1 \x \mathbb{P}^1$,  $\mathbb{P}^2$.  We have
  $\chi(S) =1,$ $\chi(\mathbb{P}^2) = 3,$ $\sigma(\mathbb{P}^2)=1, $ $\chi(\Sigma_0)= 4,$ and $\sigma(\Sigma_0) =0,$  
 so that we can consider the functions}
\begin{align*}   
\sum_{n=0}^\infty \chi(\mathrm{Hilb}^n(S))q^n = & H_{-1, 0}(q) =  1 + q + 2q^2 + 3q^3 + 5q^4 + 7q^5 +  11q^6 + \cdots \\\sum_{n=0}^\infty(-1)^n\sigma(\mathrm{Hilb}^n(\Sigma_0))q^n = & H_{0,-2}(q)  = 1 + 0q+ 2q^2 +0q^3+ 5q^4 + 0q^5+ 10q^6 + \cdots 
\\
\sum_{n=0}^\infty(-1)^n\sigma(\mathrm{Hilb}^n(\mathbb{P}^2))q^n = & H_{1,-2}(q)  = 1 - q + q^2 - 2q^3 + 3q^4 - 4q^5  + 5q^6 + \cdots .
\end{align*}



 \emph{
 Table \ref{table_0} lists   $a_2(-1, 0;n) $, $a_2(0,-2;n) $, and $a_2(1,-2;n) $ for small values of $n$,
 where  $a_N(\alpha,\beta;n)$ is the approximation obtained from Theorem \ref{thm:main_result} by summing over $1 \leq k \leq N$.  The rows correspond to the  series above in order.}

\begin{table}[h] 
\begin{center}
\begin{tabular}{|c|c|c|c|c|c|c|}\hline 
     $n$ & $1$ &$2$ & $3$  & $4$ & $5$ &$6$  \\\hline 
    $a_2({-1, 0};n)$ & $1.0029...$ & $2.0808...$ & $2.9340...$ & $5.0296...$ & $7.0278...$ & $10.9325...$ \\\hline
    $a_2({0,-2};n)$  & $0 $         & $2.1281...$ & $0$     & $4.8883...$ & $0$         & $10.1650...$ \\\hline 
    $a_2({1,-2};n)$  &  $-0.8747...$ & $1.3314...$ & $-1.9544...$ & $2.7902...$ & $-3.8958...$ & $5.3410...$  \\\hline
\end{tabular}
\end{center}
\caption{Approximate values in Theorem \ref{thm:main_result}, $N=2$} \label{table_0}
\end{table}  
\emph{As Corollary \ref{cor:asymptotics} indicates, the quality of these approximations improves as $n\rightarrow \infty$.  Moreover, choosing larger values of $N$ gives better approximations. The table below gives approximations when $N=75$. }
\begin{table}[h] 
\begin{center}
\begin{tabular}{|c|c|c|c|c|c|c|}\hline 
     $n$ & $1$ &$2$ & $3$  & $4$ & $5$ &$6$  \\\hline 
    $a_{75}({-1, 0};n)$ & $0.9999..$ & $2.0005...$ & $2.9999...$ & $4.9999...$ & $6.9999...$ & $10.9999...$ \\\hline
    $a_{75}({0,-2};n)$  & $0.0001... $ & $1.9999...$ & $-0.0002...$     & $5.0001...$ & $-0.0000...$         & $9.9999...$ \\\hline 
    $a_{75}({1,-2};n)$  &  $-1.0004...$ & $1.0003...$ & $-2.0000...$ & $3.0005...$ & $-3.9994...$ & $5.0003...$  \\\hline
\end{tabular}
\end{center}
\caption{Approximate values in Theorem \ref{thm:main_result}, $N=75$} \label{table_1}
\end{table}  
\end{example}

\begin{example} \emph{For an illustration of Corollary \ref{cor:equidistribution},  let $S'$ be $C_2 \x \bP^1$ blown up at two points, where $C_2$ is a curve of genus $2$. 
In Tables 
\ref{table_4}, and \ref{table_5}
we take }
$$\Theta^r_S(n): = \frac{b^*(r,2;n)}{\sum_{r \Mod{2}} |b^*(r,2;n)|} \qquad \text{\emph{and}} \qquad  \Theta^{r_1,r_2}_S(n): = \frac{c^*(r_1,r_2;n)}{\sum_{\substack{ r_1 \Mod{2} \\ r_2 \Mod{2}}} |c^*(r_1,r_2;n)|}.
$$
\emph{Corollary \ref{cor:equidistribution} (1) states that $\Theta^0_{S'}(n) \to 1/2$ and $\Theta^0_{S'}(n) \to 1/2$ as $n \to \infty$, which is illustrated by Table \ref{table_4}. 
Corollary \ref{cor:equidistribution} (2) states that $\Theta^{0,0}_{S'}(n),\Theta^{1,1}_{S'}(n) \to 1/4$ and $\Theta^{0,1}_{S'}(n),\Theta^{1,0}_{S'}(n) \to -1/4$, which is illustrated by Table \ref{table_5}.}

\begin{table}[h] 
\begin{center}
\begin{tabular}{|c|c|c|c|c|c|}\hline 
     $n$&$5$ & $10$ &$15$ & $20$  &$25$  \\\hline 
    $\Theta^0_{S'}(n)$  & $ 0.5714...$  &  $ 0.5054...     $  &  $ 0.4977... $  &  $ 0.4993...    $  &  $ 0.5000... $\\\hline
    $\Theta^1_{S'}(n)$& $ -0.4285... $  &  $ -0.4946... $  &  $ -0.5023...   $  &  $ -0.5006...      $  &  $ -0.5000... $ \\\hline 
\end{tabular}
\end{center}
\caption{Comparative asymptotic properties of $b^*_{S'}(r,2;n)$ } \label{table_4}
 \end{table}


\begin{table}[h] 
\begin{center}
\begin{tabular}{|c|c|c|c|c|c|}\hline 
     $n$&$5$ & $10$ &$15$ & $20$  &$25$  \\\hline 
    $\Theta^{0,0}_{S'}(n)$  & $ 0.2505... $  &  $0.2500...     $  &  $0.2500...     $  & $0.2500...     $  &$0.2500...     $  \\\hline
    $\Theta^{0,1}_{S'}(n)$&  $-0.2499... $  &  $ -0.2499... $  &  $ -0.2500...  $  & $ -0.2500...  $  &  $ -0.2500... $ \\\hline 
    $\Theta^{1,0}_{S'}(n)$&  $-0.2499... $  &  $ -0.2499... $  &  $ -0.2500...  $  & $ -0.2500...  $  &  $ -0.2500... $ \\\hline 
    $\Theta^{1,1}_{S'}(n)$& $0.2495... $  &  $  0.2499...  $  &   $  0.2499...  $  & $  0.2499...  $  & $  0.2499...  $   \\\hline
\end{tabular}
\end{center}
\caption{Comparative asymptotic properties of $c^*_{S'}(r_1,r_2;n)$ } \label{table_5}
\end{table}


\end{example}



\section*{Acknowledgements}
The authors wish to thank Ken Ono, John Duncan, and Larry Rolen for their magnificent support and
guidance. They would also like to thank Gabor Szekelyhidi for helpful communication on the topic of complex surfaces.
They also thank Emory University, 
Wesleyan University's
Marc N. Casper fund, 
the University of Notre Dame, the Asa Griggs Candler Fund, and NSF grant
DMS-1557960.

\bibliographystyle{acm}

\bibliography{main}

\end{document}